\def\updated{14 January 2005}
\count100= 1
\count101= 45

\message{: version 30nov04}

%





\font\titlefont=cmr17
\font\twelverm=cmr12
\font\ninerm=cmr9
\font\sevenrm=cmr7
\font\sixrm=cmr6
\font\fiverm=cmr5
\font\ninei=cmmi9
\font\seveni=cmmi7
\font\ninesy=cmsy9
\font\sevensy=cmsy9

\font\sixbf=cmbx6
\font\fivebf=cmbx5
\font\ninebf=cmbx9
\font\nineit=cmti9
\font\ninesl=cmsl9
\font\nineex=cmex9

\def\ninepoint{\def\rm{\fam0\ninerm}
    \textfont0 = \ninerm
    \textfont1 = \ninei
    \textfont2 = \ninesy
    \textfont3 = \nineex
    \scriptfont0 = \sevenrm
    \scriptfont1 = \seveni
    \scriptfont2 = \sevensy
    \scriptscriptfont0 = \fiverm
    \scriptscriptfont1 = \fivei
    \scriptscriptfont2 = \fivesy
    \textfont\itfam=\nineit \def\it{\fam\itfam\nineit}
    \textfont\bffam=\ninebf \scriptfont\bffam=\sixbf
    \scriptscriptfont\bffam=\fivebf \def\bf{\fam\bffam\ninebf}
    \textfont\slfam=\ninesl \def\sl{\fam\slfam\ninesl}
    \baselineskip 10pt}


\def\CC{{\rm C\kern-.18cm\vrule width.6pt height 6pt depth-.2pt
\kern.18cm}}
\def\NN{{\mathop{{\rm I}\kern-.2em{\rm N}}\nolimits}}
\def\PP{{\mathop{{\rm I}\kern-.2em{\rm P}}\nolimits}}
\def\RR{{\mathop{{\rm I}\kern-.2em{\rm R}}\nolimits}}
\def\RRt{{\titlefont I}\kern-.2em{\titlefont R}}

\def\ZZ{{\mathop{{\rm Z}\kern-.28em{\rm Z}}\nolimits}}


\def\bfm#1{{\dimen0=.01em\dimen1=.009em\makebold{$#1$}}}



\def\makebold#1{\mathord{\setbox0=\hbox{#1}%
       \copy0\kern-\wd0%
       \raise\dimen1\copy0\kern-\wd0%
       {\advance\dimen1 by \dimen1\raise\dimen1\copy0}\kern-\wd0%
       \kern\dimen0\raise\dimen1\copy0\kern-\wd0%
       {\advance\dimen1 by \dimen1\raise\dimen1\copy0}\kern-\wd0%
       \kern\dimen0\raise\dimen1\copy0\kern-\wd0%
       {\advance\dimen1 by \dimen1\raise\dimen1\copy0}\kern-\wd0%
       \kern\dimen0\raise\dimen1\copy0\kern-\wd0%
       \kern\dimen0\box0}}




%


\def\dd{\,{\rm d}} 

\def\frac#1#2{{#1 \over #2}}




\def\endproofsymbol{\makeblanksquare6{.4}}
\def\eop{\endproofsymbol\nopf}

\def\meop{~~\endproofsymbol}

\def\nopf{\medskip\goodbreak}

\def\pf{\noindent{\bf Proof: }}

\def\makeblanksquare#1#2{
\dimen0=#1pt\advance\dimen0 by -#2pt
      \vrule height#1pt width#2pt depth0pt\kern-#2pt
      \vrule height#1pt width#1pt depth-\dimen0 \kern-#1pt
      \vrule height#2pt width#1pt depth0pt \kern-#2pt
      \vrule height#1pt width#2pt depth0pt
}

\magnification\magstep0

\hsize12.1truecm\vsize18.6truecm
\hoffset.8truein



\def\title#1{\toneormore#1||||:}
\def\titexp#1#2{\hbox{{\titlefont #1} \kern-.25em%
  \raise .90ex \hbox{\twelverm #2}}\/}
\def\titsub#1#2{\hbox{{\titlefont #1} \kern-.25em%
  \lower .60ex \hbox{\twelverm #2}}\/}

\def\author#1{\bigskip\bigskip\aoneormore#1||||:\smallskip\centerline{\updated}}

\def\abstract#1{\bigskip\bigskip\medskip%
 {\ninepoint
 \narrower{\bf Abstract.~}\rm#1\bigskip
 \printtochere}\starttoc\bigskip}

\def\toneormore#1||#2||#3:{\centerline{\titlefont #1}%
    \def\next{#2}\ifx\next\empty\else\medskip\toneormore#2||#3:\fi}
\def\aoneormore#1||#2||#3:{\centerline{\twelverm #1}%
    \def\next{#2}\ifx\next\empty\else\smallskip\aoneormore#2||#3:\fi}

\newwrite\toc\def\tocone{0}\def\tochalf{.5}\def\toctwo{1}
\def\printtochere{\immediate\closeout\toc{\inputifthere{\jobname.toc}}}
\def\starttoc{\immediate\openout\toc=\jobname.toc}
\def\nexttoc#1{{\let\folio=0\edef\next{\write\toc{#1}}\next}}

\def\tocline#1#2#3{\nexttoc{\line{\hskip\parindent\noexpand\noexpand\noexpand\rm\hskip#2truecm #1\hfill#3\hskip\parindent}}}

\def\footnoterule{\kern -3pt \hrule width 0truein \kern 2.6pt}
\def\leftheadline{\ifnum\pageno=\count100 \hfill%
  \else\rm\folio\hfil\it\shortauthor\fi}
\def\rightheadline{\ifnum\pageno=\count100 \hfill%
  \else\it\shorttitle\hfil\rm\folio\fi}

\nopagenumbers
\headline{\ifodd\pageno\rightheadline \else\leftheadline\fi}
\footline{\hfil}
\null\vskip 18pt
\centerline{}
\pageno=\count100
\count102=\count100
\advance\count102 by -1
\advance\count102 by \count101


\def\copyright{\hbox{{\twelverm o}\kern-.61em\raise .46ex\hbox{\fiverm c}}}

\insert\footins{\sixrm
\medskip
\baselineskip 8pt
\leftline{Surveys in Approximation Theory
  \hfill {\ninerm \the\pageno}}
\leftline{Volume 1, 2005.
pp.~\the\pageno--\the\count102.}
\leftline{Copyright \copyright\ 2005 Surveys in Approximation Theory.}
\leftline{ISSN x-x-x}
\leftline{All rights of reproduction in any form reserved.}
\smallskip
\par\allowbreak}


\def\sect#1{\startsect\edef\showsectno{\the\sectionno}\let\tocindent\tocone%
       \soneormore#1||||:\relax\medskip\noindent\ignorespaces}

\def\soneormore#1||#2||#3:{%
   \leftline{\bf\showsectno\hskip2em #1}
   \def\next{#2}%
   \ifx\next\empty\puttocline{\showsectno\ \ #1}{\folio}%
   \else\puttocline{\showsectno\ \ #1}{}\let\showsectno\skipsectno\let\tocindent\tochalf\soneormore#2||#3:\fi}

\def\puttocline#1#2{\tocline{#1}{\tocindent}{#2}}
\def\skipsectno{\setbox0=\hbox{\the\sectionno}\hskip\wd0}

\def\subsect#1{\formal{#1.}\let\tocindent\toctwo\puttocline{#1}{\folio}}
\def\formal#1{\bigskip{\bf #1}\hskip1em}


\newcount\sectionno\sectionno0
\def\presect{\the\sectionno.}
\newcount\subsectionno
\def\startsect{\ifx\empty\presect\else\restartnums\fi%
               \subsectionno0\global\advance\sectionno by 1\relax

               \goodbreak\bigskip\smallskip}
\def\startsubsect{\global\advance\subsectionno by 1\goodbreak\bigskip}


\def\gridbox#1/#2/#3{
\vbox to #1\gridunits{#3
\ifshowgrid\tickcount=0
  \loop\cgridw%
   \vbox to 0pt{\kern\tickcount \gridunits\hrule width#2\gridunits
       height\gridwidth\vss}
   \nointerlineskip \advance\tickcount by \tickskip
   \ifdim\tickcount pt<#1pt\repeat 
  \hbox to 0pt{\tickcount=0\hbox to 0pt{\tick#1/\hss}\advance\tickcount by \tickskip%
 \loop\ifdim\tickcount pt<#2pt\nexttick\tickskip#1/\advance\tickcount by \tickskip \repeat\hss}
\else \vbox to 0pt{\hrule width#2\gridunits height0pt\vss}
\fi\vfil}\vfil}

\def\ppoint#1#2(#3,#4)#5{\setbox0=\hbox{#5}
   \dimen0=\ht0\advance\dimen0 by\dp0\divide\dimen0 by-2
   \multiply\dimen0 by#1\advance\dimen0 by#3\gridunits
   \dimen1=\wd0\divide\dimen1 by-2\multiply\dimen1 by#2
   \advance\dimen1 by#4\gridunits\dpoint(\dimen0,\dimen1){#5}}

\def\gridunits{truecm}\newcount\tickskip\tickskip1\newcount\majortick\majortick5

\def\point(#1,#2)#3{\dpoint(#1\gridunits,#2\gridunits){#3}}
\def\dpoint(#1,#2)#3{\vbox to 0pt{\kern#1
   \hbox{\kern#2{#3}}\vss}\nointerlineskip}
\newcount\rmndr   
\def\rem#1#2{\rmndr=#1{}\divide\rmndr by#2{}%
\multiply\rmndr by-#2{}\advance\rmndr by #1}
\def\cgridw{\gridwidth\finegridw{}\rem\tickcount\majortick%
   \ifnum\rmndr=0{}\gridwidth\roughgridw\fi} 

\def\tick#1/{\cgridw\vrule width\gridwidth height0pt depth#1\gridunits}
\def\nexttick#1#2/{\hbox to#1\gridunits{\hfil\tick#2/}}
\newcount\tickcount
\newdimen\finegridw\finegridw0.4pt\newdimen\roughgridw\roughgridw1.6pt
\newdimen\gridwidth
\newif\ifshowgrid \showgridtrue





\def\label#1{%
  \ifsamelabel\global\samelabelfalse\else
  \ifmmode\global\advance\eqnum by 1
  \else\global\advance\labelnum by 1
  \fi\fi
  \edef\griff{label:#1}\edef\inhalt{\lastlabel}\definieres%
  \ifmmode\eqno(\inhalt)\else\inhalt\fi
  \ifdraft\ifmmode\rlap{\fiverm #1}\else\marginal{#1}\fi\fi}

\def\eqalignlabel#1{{\def\eqno{}\let\labelnum\eqnum\label{#1}}}

\def\labelplus#1#2{\def\labelsub{#2}\relax\label{#1}\def\labelsub{}}

\def\eqalignlabelplus#1#2{{\def\eqno{}\let\labelnum\eqnum\labelplus{#1}{#2}}}


\def\recall#1{\edef\griff{label:#1}\plazieres}

\newif\ifsamelabel
\def\labelsub{}
\def\lastlabel{\presect\ifmmode\the\eqnum\else\the\labelnum\fi\labelsub}
\def\nextlabel{{\ifmmode\advance\eqnum by 1\else\advance\labelnum by 1\fi\lastlabel}}


\newcount\blackmarks\blackmarks0
\newcount\eqnum
\newcount\labelnum
\def\restartnums{\eqnum0\labelnum0}
\def\singlecount{\let\labelnum\eqnum}

 \newread\testfl
 \def\inputifthere#1{\immediate\openin\testfl=#1
    \ifeof\testfl\message{(#1 does not yet exist)}
    \else\input#1\fi\closein\testfl}

 \inputifthere{\jobname.aux}
 \newwrite\aux
 \immediate\openout\aux=\jobname.aux

\def\plazieres{\expandafter\ifx\csname\griff\endcsname\relax%
  \xdef\esfehlt{\griff}\blackmark\else{\csname\griff\endcsname}\fi}

\def\definieres{\expandafter\xdef\csname\griff\endcsname{\inhalt}%
 \def\blankkk{ }\expandafter\immediate\write\aux{%
 \string\expandafter\def\string\csname%
 \blankkk\griff\string\endcsname{\inhalt}}}

\def\blackmark{\ifnum\blackmarks=0\global\blackmarks=1%
 \write16{============================================================}%
 \write16{Some forward reference is not yet defined. Re-TeX this file!}%
 \write16{============================================================}%
 \fi\immediate\write16{undefined forward reference: \esfehlt}%
 {\vrule height10pt width2pt depth2pt}\esfehlt%
 {\vrule height10pt width2pt depth2pt}}

\def\marginal#1{\strut\vadjust{\kern-\strutdepth%
\vtop to \strutdepth{\baselineskip\strutdepth\vss\llap{\fiverm#1\ }\null}}}
\def\strutdepth{\dp\strutbox}


\newif\ifdraft

\newcount\hour\newcount\minutes
\def\draft{\drafttrue
\def\comment##1{{\bf comment: ##1}}
\headline={\sevenrm \hfill\ifx\filenamed\undefined\jobname\else\filenamed\fi%
(.tex) (as of \ifx\updated\undefined???\else\updated\fi)
 \TeX'ed at {\hour\time\divide\hour by 60{}%
\minutes\hour\multiply\minutes by 60{}%
\advance\time by -\minutes
\the\hour:\ifnum\time<10{}0\fi\the\time\  on \today\hfill}}
}

\def\today{\number\day\space%
\ifcase\month\or January\or February\or March\or April\or May\or June\or
 July\or August\or September\or October\or November\or December\fi%
\space\number\year}



\def\References{\goodbreak\bigskip\centerline{\bf References}%
   \tocline{\skipsectno\ \  References}{\tocone}{\folio}%
   \bigskip\frenchspacing}

\def\bibitem{\smallskip\noindent}


\gdef\formfirstauthor{\the\lastname, \the\firstname}
\gdef\formnextauthor{, \the\firstname\the\lastname}
\gdef\formotherauthor{ and \the\firstname\the\lastname}
\gdef\formlastauthor{,\formotherauthor}

\gdef\formB{\the\au\ [\yr] ``\the\ti'', \the\pb, \the\pl. \setcitelabel}
\gdef\formD{\the\au\ [\yr] ``\the\ti'', dissertation, \the\pl. \setcitelabel}
\gdef\formJ{\the\au\ [\yr] \the\ti, {\sl\the\jr}\ifx\vl\empty%
\else\ {\bf\vl}\fi, \pp. \setcitelabel}
\gdef\formP{\the\au\ [\yr] \the\ti, in {\sl\the\tit},
\getfirstchar\aut\ifx\firstchar\unknownx\else\the\aut, ed\edsop, \fi
\getfirstchar\pub\ifx\firstchar\unknownx\else\the\pub, \fi \the\pl, \pp. \setcitelabel}
\gdef\formR{\the\au\ [\yr] \the\ti\ifx\is\empty\else, \is\fi. \setcitelabel}

\newtoks\lastname
\newtoks\firstname
\newtoks\au
\newtoks\aut
\newtoks\ti
\newtoks\tit
\newtoks\pb
\newtoks\pub
\newtoks\pl
\newtoks\jr

\newtoks\rhlau

\def\setcitelabel{\edef\griff{cit\rh}\edef\inhalt{\the\rhlau\ \yr}\definieres}
\def\setcitelabel{}

\def\getfirstchar#1{\edef\theword{\the#1}\expandafter\getit\theword:}
\def\getit#1#2:{\def\firstchar{#1}}
\def\unknownx{x}

\newif\ifonesofar
\def\concat#1{\edef\audef{{#1}}\au=\audef}
\def\decodeauthor#1, #2,#3;{\lastname={#1}\firstname={#2}%
\concat{\formfirstauthor}\onesofartrue%
\def\morerhlau{}%
\def\next{#3}\ifx\next\empty\else\def\morerhlau{ et al.}\decodemoreauthor#3;\fi
\edef\morerhlauu{{\the\lastname\morerhlau}}\rhlau=\morerhlauu}
\def\decodemoreauthor#1, #2,#3;{\lastname={#1}\firstname={#2}%
\def\next{#3}\ifx\next\empty\let\formaut=\formlastauthor%
\ifonesofar\ifx\formotherauthor\undefined\else\let\formaut=\formotherauthor%
\fi\fi\concat{\the\au\formaut}%
\else\onesofarfalse\concat{\the\au\formnextauthor}\decodemoreauthor#3;\fi}

\def\refproc #1(#2; #3; {\decodeproc#2; \xdef\yr{#3}}
\def\decodeproc#1), #2 (ed#3.),#4 (#5); {%
 \global\tit={#1}\global\aut={#2}\xdef\edsop{#3}\global
 \pub={#4}\global\pl={#5}}


\def\refB #1; #2; #3 (#4); #5; {\decodeauthor#1,;%
   \ti={#2}\pb={#3}\pl={#4}\def\yr{#5}\bibitem\formB}

\def\refD #1; #2; #3; #4; {\decodeauthor#1,;%
   \ti={#2}\pl={#3}\def\yr{#4}\bibitem\formD}

\def\refJ #1; #2; #3; #4; #5; #6; {\decodeauthor#1,;%
    \ti={#2}\jr={#3}\def\vl{#4}\def\yr{#5}\def\pp{#6}\bibitem\formJ}

\def\refP #1; #2; #3; #4; {{\global\aut={\vrule height15pt width15pt depth0pt}%
 \global\tit={{\bf the specified proceedings does not exist in our files}}%
 \xdef\edsop{}\global\pub={}\def#3{}\input proceed }\decodeauthor#1,;%
        \ti={#2}\def\pp{#4}\bibitem\formP}

\def\refQ #1; #2; (#3; #4; #5; {\decodeproc#3; \decodeauthor#1,;%
   \ti={#2}\def\yr{#4}\def\pp{#5}\bibitem\formP}

\def\refR #1; #2; #3; #4; {\decodeauthor#1,;%
         \ti={#2}\def\is{#3}\def\yr{#4}\bibitem\formR}




\title{Density in Approximation Theory}
\author{Allan Pinkus}

\def\shorttitle{Density in Approximation Theory}
\def\shortauthor{Allan Pinkus}

\def\alp{\alpha}                
\def\gam{\gamma}                
\def\del{\delta}                
\def\eps{\varepsilon}
\def\lam{\lambda}               \def\Lam{\Lambda}
\def\sig{\sigma}                
\def\ome{\omega}                \def\Ome{\Omega}

\def\calA{{\cal A}}
\def\calM{{\cal M}}
\def\calN{{\cal N}}
\def\calP{{\cal P}}
\def\calS{{\cal S}}

\def\bfa{{\bf a}}               
\def\bfm{{\bf m}}               
\def\bfw{{\bf w}}               
\def\bfx{{\bf x}}               
\def\bfy{{\bf y}}               

\def\CC{{\rlap {\raise 0.4ex \hbox{$\scriptscriptstyle |$}}\hskip -0.13em C}}

\def\NN{{I\!\!N}}
\def\PP{{I\hskip-2pt P}}
\def\QQ{{\rlap {\raise 0.4ex \hbox{$\scriptscriptstyle |$}}\hskip -0.1em Q}}
\def\RR{{I\!\!R}}
\def\ZZ{{Z\!\!\! Z}}

\def\incl{\subseteq}

\def\nek{,\ldots,}

\def\sdp{\times \hskip -0.3em {\raise 0.3ex\hbox{$\scriptscriptstyle |$}}} 
\def\union{\bigcup}
\def\amp{\raise 6pt\hbox{$\scriptstyle\bf1$}}
\def\scirc{{\, \raise 0.3ex\hbox{$\scriptstyle \circ$}\,}}
\def\divdif{\mathord{\kern.43em\vrule width.6pt height4.95pt
depth.28pt\kern-.43em\Delta}}

\def\span{{\rm span}}
\def\supp{{\rm supp}}

\def\\{{\backslash}}

\def\oA{{\overline A}}

\def\hatg{{\widehat g}}
\def\hatphi{{\widehat\phi}}
\def\hatmu{{\widehat\mu}}
\def\hatpsi{{\widehat\psi}}

\def\tilC{{\widetilde C}}

\overfullrule=0pt


\abstract {Approximation theory is concerned with the ability to
approximate functions by simpler and more easily calculated
functions. The first question we ask in approximation theory
concerns the {\it possibility of approximation}. Is the given
family of functions from which we plan to approximate dense in
the set of functions we wish to approximate? In this work we
survey some of the main density results and density methods.}


\sect{Introduction}
Approximation theory is that area of analysis which, at its core,
is concerned with the ability to approximate functions by simpler
and more easily calculated functions. It is an area which, like
many other fields of analysis, has its primary roots in the
mathematics of the 19th century.

At the beginning of the 19th century functions were essentially
viewed via concrete formulae, series, or as solutions of
equations. However largely as a consequence of the claims of
Fourier and the results of Dirichlet, the modern concept of a
function distinguished by its requisite properties was introduced
and accepted. Once a function, and more specifically a continuous
function, is defined implicitly rather than explicitly, the birth
of approximation theory becomes an inevitable and unavoidable
development.

It is in the theory of Fourier series that we find some of the
first results of approximation theory. These include conditions on
a function that ensure the pointwise or uniform convergence (of
the partial sums) of its Fourier series, as well as the
omnipresent $L^2$-convergence. Similar results were also developed
for other orthogonal series, and for power series (analytic
functions). However these results are of a rather particular form.
They are concerned with conditions for when certain formulae hold.
In the classical theory of Fourier series one does not ask if
trigonometric polynomials can be used to approximate, or even if
the information provided by the Fourier coefficients is sufficient
to provide an approximation. Rather one wants to know if the
partial sums of the Fourier series converge to the function in
question.

The first question we ask in approximation theory concerns the
{\it possibility of approximation}. Is the given family of
functions from which we plan to approximate dense in the set of
functions we wish to approximate? That is, can we approximate any
function in our set, as well as we might wish, using arbitrary
functions from our given family? In this work we survey some of
the main density results and density methods.

The first significant density results were those of Weierstrass
who proved in 1885 (when he was 70 years old!) the density of
algebraic polynomials in the class of continuous real-valued
functions on a compact interval, and the density of trigonometric
polynomials in the class of $2\pi$-periodic continuous real-valued
functions. These theorems were, in a sense, a counterbalance to
Weierstrass' famous example of 1861 on the existence of a
continuous nowhere differentiable function. The existence of such
functions accentuated the need for analytic rigour in mathematics,
for a further understanding of the nature of the set of continuous
functions, and substantially influenced the further development of
analysis. If this example represented for some a `lamentable
plague' (as Hermite wrote to Stieltjes on May 20, 1893, see
Baillard and Bourget [1905]), then the approximation theorems were
a panacea. While on the one hand the set of continuous functions
contains deficient functions, on the other hand every continuous
function can be approximated arbitrarily well by the ultimate in
smooth functions, the polynomials.

The Weierstrass approximation theorems spawned numerous
generalizations which were applied to other families of functions.
They also led to the development of two general methods for
determining density. These are the Stone-Weierstrass theorem
generalizing the Weierstrass theorem to subalgebras of $C(X)$, $X$
a compact space, and the Bohman-Korovkin theorem characterizing
sequences of positive linear operators that approximate the
identity operator, based on easily checked, simple, criteria.

A different and more modern approach to density theorems is via
``soft analysis''. This functional analytic approach actually
dates back almost 100 years. A linear subspace $M$ of a normed
linear space $E$ is dense in $E$ if and only if the only
continuous linear functional that vanishes on $M$ is the
identically zero functional. For the space $C[a,b]$ this result
can already be found in the work of F.~Riesz from 1910 and 1911 as
one of the first applications of his ``representation theorem''
characterizing the set of all continuous linear functionals on
$C[a,b]$.

Density theorems can be found almost everywhere in analysis, and
not only in analysis. (For a density result equivalent to the
Riemann Hypothesis see Conrey [2003, p.~345].) In this article we
survey some of the main results regarding density of linear
subspaces in spaces of continuous real-valued functions endowed
with the uniform norm. We only present a limited sampling of the
many, many density results to be found in approximation theory and
in other areas. A monograph many times the length of this work
would not suffice to include all results. In addition, we do not
prove all the results we quote. Writing a paper such as this
involves compromises. We hope, nonetheless, that you the reader
will find something here of interest.

\sect{The Weierstrass Approximation Theorems} We first fix some
notation. We let $C[a,b]$ denote the class of continuous
real-valued functions on the closed interval $[a,b]$, and
$\tilC[a,b]$ the class of functions in $C[a,b]$ satisfying
$f(a)=f(b)$. ($\tilC[a,b]$ may be regarded as the restriction to
$[a,b]$ of $(b-a)$-periodic functions in $C(\RR)$.) We
denote by $\Pi_n$ the space of algebraic polynomials of degree at
most $n$, i.e., $$\Pi_n =\span\{ 1,x\nek x^n\},$$ and by $T_n$ the
space of trigonometric polynomials of degree at most $n$, i.e.,
$$T_n =\span \{1, \sin x, \cos x\nek \sin nx, \cos nx \}.$$

The paper stating and proving what we call the Weierstrass
approximation theorems is Weierstrass [1885]. It seems that the
importance of the paper was immediately appreciated, as the paper
appeared in translation (in French) one year later in Weierstrass
[1886]. Weierstrass was interested in complex function theory and
in the ability to represent functions by power series and function
series. He viewed the  results obtained in this 1885 from that
perspective. The title of the paper emphasizes this viewpoint. The
paper is titled {\it On the possibility of giving an analytic
representation to an arbitrary function of a real variable}. We
state the Weierstrass theorems, not as given in his paper, but as
they are currently stated and understood.

\proclaim Weierstrass Theorem \label{weier1}.  For every finite
$a<b$ algebraic polynomials are dense in $C[a,b]$. That is, given
an $f$ in $C[a,b]$ and an arbitrary $\eps>0$ there exists an
algebraic polynomial $p$ such that $$|f(x)-p(x)|<\eps$$ for all
$x$ in $[a,b]$.\nopf

\proclaim Weierstrass Theorem \label{weier2}. Trigonometric
polynomials are dense in $\tilC[0,2\pi]$. That is, given an $f$ in
$\tilC[0,2\pi]$ and an arbitrary $\eps>0$ there exists a
trigonometric polynomial $t$ such that $$|f(x)-t(x)|<\eps$$ for
all $x$ in $[0,2\pi]$.\nopf

These are the first significant density theorems in analysis. They
are generally paired since in fact they are equivalent. That is,
each of these theorems follows from the other.

It is interesting to read this paper of Weierstrass, as his perception
of these approximation theorems was most certainly different from
ours. Weierstrass' view of analytic functions was of functions
that could be represented by power series. The approximation
theorem, for him, was an extension of this result to continuous
functions.

Explicitly, let $(\eps_n)$ be any sequence of positive values for
which $\sum_{n=1}^\infty \eps_n <\infty$. Let $p_n$ be an
algebraic polynomial (which exists by Theorem \recall{weier1})
satisfying $$ \|f-p_n\| :=\max_{a\le x\le b}
|f(x)-p_n(x)|<\eps_n,\qquad n=1,2,...$$ Set $q_0=p_1$ and $q_n=
p_{n+1}-p_n$, $n=1,2,...$ Then
$$f(x) = \sum_{n=0}^\infty q_n(x).$$ Thus every continuous
function can be represented by a polynomial series that converges
both absolutely and uniformly. Similarly, `nice' functions in
$\tilC[0,2\pi]$ enjoy the property that their Fourier series
converges absolutely and uniformly. What Weierstrass proved
was that every function in $\tilC[0,2\pi]$ can be represented by a
trigonometric polynomial series that converged both absolutely and
uniformly.

The paper Weierstrass [1885] was reprinted in Weierstrass'
Mathematische Werke (collected works) with some notable additions.
While this reprint appeared in 1903, there is reason to assume
that Weierstrass himself edited this paper. One of these additions
was a short ``introduction''. We quote it (verbatim in meaning if
not in fact).

{\it The main result of this paper, restricted to the one variable
case, can be summarized as follows:

Let $f\in C(\RR)$. Then there exists a sequence $f_1, f_2, \ldots$
of entire functions for which $$f(x)=\sum_{i=1}^\infty f_i(x)$$
for each $x\in \RR$. In addition the convergence of above sum is
uniform on every finite interval.}

Note that there is no mention of the fact that the $f_i$ may be
assumed to be polynomials.

Weierstrass' proof of Theorem \recall{weier1}\ is rather
straightforward. The same is not quite true of his proof of
Theorem \recall{weier2}. He extends $f$ from $[a,b]$ so that it is
continuous and bounded on all of $\RR$. He then smooths $f$ by
convolving it with the normalized heat (Gauss) kernel
$(1/k\sqrt{\pi})e^{-(x/k)^2}$. This ``smoothed'' $f_k$ is entire
and is therefore uniformly approximable on the finite interval
$[a,b]$ by its truncated power series. Moreover the $f_k$
uniformly approximate $f$ on $[a,b]$ as $k\to 0+$. Together this
implies the desired result.

Over the next twenty-five or so years numerous alternative proofs
were given to one or the other of these two Weierstrass results by
a roster of some of the best analysts of the period. The proofs
use diverse ideas and techniques. There are the proofs by
Weierstrass, Picard, Fej\'er, Landau and de la Valle\'e Poussin
that used singular integrals, proofs based on the idea of
approximating one particular function by Runge (Phragm\'en),
Lebesgue, Mittag-Leffler, and Lerch, proofs based on Fourier
series by Lerch, Volterra and Fej\'er, and the wonderful proof of
Bernstein. Details concerning all these proofs can be found, for
example, in Pinkus [2000] and Pinkus [2005]. We explain, without
going into all the details, three of these proofs.

One of the more elegant and cited proofs of Weierstrass' theorem
is due to Lebesgue [1898]. This was Lebesgue's first published
paper. He was, at the time of publication, a 23 year old student
at the \'Ecole Normale Sup\'erieure. The idea of his proof is
simple and useful. Lebesgue noted that each $f$ in $C[a,b]$ can be
easily approximated by a continuous, piecewise linear curve
(polygonal line). Each such polygonal line is a linear combination
of translates of $|x|$. As algebraic polynomials (of any fixed
degree) are translation invariant, it thus suffices to prove that
one can uniformly approximate $|x|$ arbitrarily well by
polynomials on any interval containing the origin. Lebesgue then
does exactly that. Explicitly $$|x| = 1 -\sum_{n=1}^\infty
a_n(1-x^2)^n$$ where $a_1= 1/2$, and $$a_n = {{(2n-3)!}\over
{2^{2n-2} n! (n-1)!}},\qquad n=2,3,\ldots$$ This ``power series''
converges absolutely and uniformly to $|x|$ for all $|x|\le 1$.
Truncating this series we obtain a series of polynomial
approximants to $|x|$.

When Fej\'er was 20 years old he published Fej\'er [1900] that
formed the basis for his doctoral thesis. Fej\'er proved more than
the Weierstrass approximation theorem (for trigonometric
polynomials). He proved that for any $f$ in $\tilC[0,2\pi]$ it is
possible to uniformly approximate $f$ based solely on the
knowledge of its Fourier coefficients. He did not obtain this
approximation by taking the partial sums of the Fourier series. It
is well-known that these do not necessarily converge. Rather, he
obtained it by taking the Ces\`aro sums of the partial sums of the
Fourier series. In other words, assume that we are given the
Fourier series of $f$ $$f(x) \sim \sum_{k=-\infty}^\infty c_k
e^{ikx},$$ where $$c_k={1\over {2\pi}}\int_{0}^{2\pi}
f(x)e^{-ikx}\dd x$$ for every $k\in \ZZ$. Define the $n$th partial
sums of the Fourier series via $$s_n(x) := \sum_{k=-n}^n c_k
e^{ikx}$$ and set $$\sig_n(f;x) ={{s_0(x)+\cdots + s_n(x)}\over
{n+1}}.$$ The $\sig_n$ are termed the $n$th Fej\'er operator. Note
that $\sig_n(f;\cdot)$ belongs to $T_n$ for each $n$. What Fej\'er
proved was that, for each $f$ in $\tilC[0, 2\pi]$,
$\sig_n(f;\cdot)$ tends uniformly to $f$ as $n\to \infty$. This
was also the first proof which used a specifically given sequence
of {\it linear} operators.

Simpler linear operators that approximate were introduced by
Bernstein [1912/13]. These are the Bernstein polynomials. For $f$
in $C[0,1]$ they are defined by $$B_n(f\,;\,x) = \sum_{m=0}^n
f\left({m \over n}\right){m \choose n}x^m(1-x)^{n-m}.$$ Bernstein
proved, by probabilistic methods, that the $B_n(f\,;\cdot)$
converge uniformly to $f$ as $n \to\infty$. A proof of this
convergence is to be found in Example 4.2.

\sect{The Functional Analytic Approach}
The Riesz representation theorem characterizing the space of
continuous linear functionals on $C[a,b]$ is contained in the 1909
paper of F.~Riesz [1909]. The following year, in a rarely
referenced paper, Riesz [1910] also announced the following
(stated in more modern terminology).

\proclaim Theorem \label{riesz1}. Let $u_k\in C[a,b]$, $k\in K$,
where $K$ is an index set. A necessary and sufficient condition
for the existence of a continuous linear functional $F$ on
$C[a,b]$ satisfying $$F(u_k) = c_k, \qquad k\in K$$ with $\|F\|\le
L$ is that $$\left|\sum_{k\in K'} a_k c_k \right| \le L
\left\|\sum_{k\in K'} a_k u_k\right\|_{\infty}$$ hold for every
finite subset $K'$ of $K$, and all real $a_k$. \nopf

In this same paper Riesz also states the parallel result for
$L^p[a,b]$, $1< p <\infty$. Questions concerning existence and
uniqueness in moment problems were of major importance in the
development of functional analysis. The full details of the 1910
announcement appear in Riesz [1911]. In these papers is also to be
found the following result (again we switch to more modern
terminology).

\proclaim Theorem \label{riesz2}. Let $M$ be a linear subspace of
$C[a,b]$. Then $f\in C[a,b]$ is in the closure of $M$, i.e., $f$
can be uniformly approximated by elements of $M$, if and only if
every continuous linear functional on $C[a,b]$ that vanishes on
$M$ also vanishes on $f$.\nopf

Riesz quotes E.~Schmidt as the author of the {\it very interesting
problem} whose solution is the above Theorem \recall{riesz2}. As
he writes, the question asked is: {\it Being given a countable
system of functions $\phi_n\in C[a,b]$, $n=1, 2,...$, how can one
know if one can approximate arbitrarily and uniformly every $f\in
C[a,b]$ by the $\phi_n$ and their linear combinations?} (Riesz
[1911, p.~51]). Schmidt, in his thesis in Schmidt [1905], had
given both a necessary and a sufficient condition for the above to
hold. Both were orthogonality type conditions. However neither was
the correct condition. The concept of a linear functional
vanishing on a set of functions is very orthogonal-like. Lerch's
theorem (Lerch [1892], see also the more accessible Lerch [1903]),
states that if $h\in C[0,1]$ and $$\int_0^1 x^n h(x)\dd x=0,\qquad
n=0,1\nek$$ then $h=0$. This theorem was well-known and frequently
quoted. So it was not unreasonable to look for conditions of the
form given in Theorem \recall{riesz2}.

Lerch's theorem is, in fact, a simple consequence of Weierstrass'
theorem. If $p_k$ is a sequence of polynomials that uniformly
approximate $h$, then $$\lim_{k\to\infty} \int_0^1 p_k(x)h(x)\dd x=
\int_0^1 [h(x)]^2 \dd x.$$ However $$ \int_0^1 p_k(x)h(x) \dd x=0$$
for every $k$, and thus $$\int_0^1 [h(x)]^2 \dd x =0$$ which, since
$h$ is continuous, implies $h=0$.

As Riesz states, one consequence of the above Theorem
\recall{riesz2} is that $M$ is dense in $C[a,b]$ if and only if no
nontrivial continuous linear functional vanishes on $M$. The proof
of Theorem \recall{riesz2}, contained in Riesz [1911], is just an
application of Theorem \recall{riesz1}.

\medskip
\pf  We start with the simple direction. Assume $f$ is in
the closure of $M$. If $F$ is a continuous linear functional that
vanishes on $M$, then $$F(f) = F(f-g)$$ for every $g\in M$. Given
$\eps>0$, there exists a $g^*\in M$ for which $$\|f-g^*\|_\infty
<\eps.$$ Thus $$|F(f)|=|F(f-g^*)|\le \|F\| \|f-g^*\|_\infty <
\eps \|F\|.$$ As this is valid for every $\eps>0$ we have
$F(f)=0$.

Now assume that $f$ is not in  the closure of $M$. Thus
$\|f-g\|_{\infty}\ge d > 0$ for every $g\in M$. From this
inequality and Theorem \recall{riesz1} there necessarily exists a
continuous linear functional $F$ on $C[a,b]$ satisfying $F(g) =
0$, for all $g\in M$, $F(f) =1$, and  $\|F\|\le L$ for any $L \ge
1/d$. This holds since we have $$|a| \le Ld |a|\le L
\|af-g\|_{\infty},$$ for all $g \in M$ and all $a$. \eop

Shortly thereafter Helly [1912] applied these results to a
question concerning the range of an integral operator. He proved
the following two theorems.

\proclaim Theorem \label{}.  Let $K\in C([a,b]\times [a,b])$ and
$f\in C[a,b]$. Then a necessary and sufficient condition for the
existence of a measure of bounded total variation $\nu$ satisfying
$$f(x) = \int_a^b K(x,y) \dd \nu(y),$$ is the existence of a
constant $L$ for which $$\left| \sum_{k=1}^n a_k f(x_k)\right| \le
L \left| \sum_{k=1}^n a_kK(x_k, y) \right|$$ for all points
$x_1\nek x_n$ in $[a,b]$, all real values $a_1\nek a_n$, all $y\in
[a,b]$, and all $n$.\nopf

\proclaim Theorem \label{}.  Let $K\in C([a,b]\times [a,b])$.
Then a necessary and sufficient condition for an $f\in C[a,b]$ to
be uniformly approximated by functions of the form $$ \int_a^b
K(x,y) \phi(y) \dd y$$ where the $\phi$ are piecewise continuous
functions, is that for every measure $\mu$ of bounded total
variation satisfying $$\int_a^b K(x,y) \dd \mu(x) =0$$ we also have
$$\int_a^b f(x) \dd \mu(y) =0.$$\nopf

In 1911 the concept of a normed linear space did not exist, and
the Hahn-Banach theorem had yet to be discovered (although the
Helly [1912] paper contains results that come close). Banach's
proof of the Hahn-Banach theorem appears in Banach [1929] (Hahn's
appears in Hahn [1927]). Both the Hahn and Banach papers contain a
general form of Theorems \recall{riesz1}, namely the Hahn-Banach
theorem. Both also essentially contain the statement that a linear
subspace is dense in a normed linear space if and only if no
nontrivial continuous linear functional vanishes on the subspace.
Banach, in his book Banach [1932, p.~57], prefaces these next two
theorems with the statement: {\it We are now going to establish
some theorems that play in the theory of normed spaces the
analogous role to that which the Weierstrass theorem on the
approximation of continuous functions by polynomials plays in the
theory of functions of a real variable.}

\proclaim Theorem \label{banach1}.  Let $M$ be a linear subspace
of a real normed linear space $E$. Assume $f\in E$ and
$$\|f-g\|\ge d >0$$ for all $g\in M$. Then there exists a
continuous linear functional $F$ on $E$ such that $F(g)=0$ for all
$g\in M$, $F(f)=1$, and $\|F\| \le 1/d$.\nopf

The result of Theorem \recall{banach1} replaces Theorem
\recall{riesz1} in the proof of Theorem \recall{riesz2} to give us
the well-known

\proclaim Theorem \label{banach2}. Let $M$ be a linear subspace of
a real normed linear space $E$. Then $f\in E$ is in the closure of
$M$ if and only if every continuous linear functional on $E$ that
vanishes on $M$ also vanishes on $f$.

In none of these works of Hahn and Banach are the above-mentioned
1910 or 1911 papers of Riesz mentioned. These Riesz papers seem to
have been essentially forgotten. In fact the general method of
proof of density based on this approach is to be found in the
literature only after the appearance of the book of Banach and the
blooming of functional analysis. The name of Riesz is often
mentioned in connection with this method, but only because of the
Riesz representation theorem and similar duality results bearing
his name.

Today we also recognize the Hahn-Banach theorem as a separation
theorem, and as such we also have the following two results.

\proclaim Theorem \label{3.7}.  Let E be a real normed linear
space, $\phi_n$ elements of $E$, $n\in I$, and $f\in E$. Then $f$
may be approximated by finite {\it convex} linear combinations of
the $\phi_n$, $n\in I$, if and only if $$\sup\{F(\phi_n)\,:\,n\in
I\} \ge F(f)$$ for every continuous linear functional (form) $F$
on $E$.\nopf

\proclaim Theorem \label{3.8}. Let E be a real normed linear
space, $\phi_n$ elements of $E$, $n\in I$, and $f\in E$. Then $f$
may be approximated by finite {\it positive} linear combinations
of the $\phi_n$, $n\in I$, if and only if for every  continuous
linear functional (form) $F$ on $E$ satisfying $F(\phi_n)\ge 0$
for every $n\in I$ we have $F(f)\ge 0$.\nopf

Theorem \recall{3.8} follows from Theorem \recall{3.7} by
considering the convex cone generated by the $\phi_n$.

There are numerous generalizations of these results. The book of
Nachbin [1967] where these results may be found is one of the few
to concentrate on density theorems. Much of the book is taken up
with the Stone-Weierstrass theorem. However there are also other
results such as the above Theorems \recall{3.7} and \recall{3.8}.

\sect{Other Density Methods}
The Weierstrass theorems had a significant influence on the
development of density results, even though the theorems
themselves simply prove the density of algebraic and trigonometric
polynomials in the appropriate spaces. Various proofs of the
Weierstrass theorems, for example, provided insights that led to
the development of two general methods for determining density. We
discuss these methods in this section.

The first of these methods is given by the Stone-Weierstrass
theorem. This theorem was originally proven in Stone [1937]. Stone
subsequently reworked his proof in Stone [1948]. It represents, as
stated by Buck [1962, p.~4], {\it one of the first and most
striking examples of the success of the algebraic approach to
analysis}. There have since been numerous modifications and
extensions. See, for example, Nachbin [1967], Prolla [1993] and
references therein.

We recall that an {\it algebra} is a linear space on which
multiplication between elements has been suitably defined
satisfying the usual commutative and associative type postulates.
Algebraic and trigonometric polynomials in any finite number of
variables are algebras. A set in $C(X)$ {\it separates points} if
for any distinct points $x,y\in X$ there exists a $g$ in the set
for which $g(x)\ne g(y)$.

\proclaim Stone-Weierstrass Theorem \label{}. Let $X$ be a compact
set and let $C(X)$ denote the space of continuous real-valued
functions defined on $X$. Assume $A$ is a subalgebra of $C(X)$.
Then $A$ is dense in $C(X)$ in the uniform norm if and only if $A$
separates points and for each $x\in X$ there exists an $f\in A$
satisfying $f(x)\ne 0$.

\pf The necessity of the two conditions is obvious. We
prove the sufficiency.

First some preliminaries. From the Weierstrass theorem we have the
existence of a sequence of algebraic polynomials $(p_n)$ (with
constant term zero) that uniformly approximates the function $|t|$
on $[-c,c]$, any $c>0$. As such, if $f$ is in $\oA$, the closure
of $A$ in the uniform norm, then so is $p_n(f)$ for each $n$ which
implies that $|f|$ is also in $\oA$. Furthermore $$\max\{ f(x),
g(x) \} = {{f(x)+g(x) + |f(x)-g(x)|}\over 2}$$ and $$\min\{ f(x),
g(x) \} = {{f(x)+g(x) - |f(x)-g(x)|}\over 2}.$$ It thus follows
that if $f,g\in \oA$, then $\max\{f,g\}$ and $\min\{f,g\}$ are
also in $\oA$. This of course extends to the maximum and minimum
of any finite number of functions.

Finally, let $x, y$ be any distinct points in $X$, and
$\alp,\beta\in \RR$. We claim that there exists an $h\in A$
satisfying the interpolation conditions $h(x)=\alp$ and
$h(y)=\beta$. By assumption there exists a $g\in A$ for which
$g(x)\ne g(y)$, and functions $f_1$ and $f_2$ in $A$ such that
$f_1(x)\ne 0$ while $f_2(y)\ne 0$. If $g(x)=0$ then we can
construct the desired $h$ as a linear combination of $g$ and
$f_1$. Similarly, if $g(y)=0$ then we can construct the desired
$h$ as a linear combination of $g$ and $f_2$. Assuming $g(x)$ and
$g(y)$ are both not zero, the desired $h$ can be constructed, for
example, as a linear combination of $g$ and $g^2$.

We now present a proof of the theorem. Given $f\in C(X)$, $\eps>0$
and $x\in X$, for every $y\in X$ let $h_y\in A$ satisfy
$h_y(x)=f(x)$ and $h_y(y)=f(y)$. Since $f$ and $h_y$ are
continuous there exists a neighborhood $V_y$ of $y$ for which
$h_y(w) \ge f(w) -\eps$ for all $w\in V_y$. The $\cup_{y\in X}
V_y$ cover $X$. As $X$ is compact, it has a finite subcover, i.e.,
there are points $y_1\nek y_n$ in $X$ such that $$\union_{i=1}^n
V_{y_i}=X.$$ Let $g=\max\{h_{y_1}\nek h_{y_n}\}$. Then $g\in
\oA$ and $g(w) \ge f(w)-\eps$ for all $w\in X$.

The above $g$ depends upon $x$, so we shall now denote it by
$g_x$. It satisfies $g_x(x)=f(x)$ and $g_x(w)\ge f(w)-\eps$ for
all $w\in X$. As $f$ and $g_x$ are continuous there exists a
neighborhood $U_x$ of $x$ for which $g_x(w) \le f(w) +\eps$ for
all $w\in U_x$. Since $\cup_{x\in X} U_x$ covers $X$, it has a
finite subcover. Thus there exist points $x_1\nek x_m$ in $X$ for
which $$\union_{i=1}^m U_{x_i}=X.$$

Let $$F=\min\{g_{x_1}\nek g_{x_m}\}.$$ Then $F\in \oA$ and
$$f(w)-\eps\le F(w) \le f(w)+\eps$$ for all $w\in X$. Thus
$$\|f-F\|\le \eps.$$ This implies that $f\in \oA$. \eop

\noindent {\bf Example 4.1.} As we mentioned prior to the
statement of the Stone-Weierstrass theorem, algebraic polynomials
in any finite number of variables form an algebra. They also
separate points and contain the constant function. Thus algebraic
polynomials in $m$ variables are dense in $C(X)$ where $X$ is any
compact set in $\RR^m$. This fact first appeared in print (at
least for squares) in Picard [1891] which also contains an
alternative proof of Weierstrass' theorems. The paper Weierstrass
[1885] as ``reprinted'' in Weierstrass' Mathematische Werke in
1903 contains an additional 10 pages of material including a proof
of this multivariable analogue of his theorem.

\medskip
Another method that can be used to prove density is based on what
is called the Korovkin theorem or the Bohman-Korovkin theorem. A
primitive form of this theorem was proved by Bohman in Bohman
[1952]. His proof, and the main idea in his approach, was a
generalization of Bernstein's proof of the Weierstrass theorem.
Korovkin one year later in Korovkin [1953] proved the same theorem
for integral type operators. Korovkin's original proof is in fact
based on positive singular integrals and there are very obvious
links to Lebesgue's work on singular operators that, in turn, was
motivated by various of the proofs of the Weierstrass theorems.
Korovkin was probably unaware of Bohman's result. Korovkin
subsequently much extended his theory, major portions of which can
be found in his book Korovkin [1960]. The theorem and proof as
presented here is taken from Korovkin's book.

A linear operator $L$ is {\it positive} (monotone) if $f\ge 0$
implies $L(f)\ge 0$.

\proclaim Bohman--Korovkin Theorem \label{}.  Let $(L_n)$ be a
sequence of positive linear operators mapping $C[a,b]$ into
itself. Assume that $$\lim_{n\to\infty} L_n(x^i) = x^i,\qquad
i=0,1,2,$$ and the convergence is uniform on $[a,b]$. Then
$$\lim_{n\to\infty} (L_nf)(x) = f(x)$$ uniformly on $[a,b]$ for
every $f\in C[a,b]$.

\pf Let $f\in C[a,b]$. As $f$ is uniformly continuous,
given $\eps>0$ there exists a $\del>0$ such that if
$|x_1-x_2|<\del$ then $|f(x_1)-f(x_2)|< \eps$.

For each $y\in [a,b]$, set $$p_u(x) = f(y) +\eps +{{2\|f\|
(x-y)^2}\over {\del^2}}$$ and $$p_\ell(x) = f(y) -\eps -{{2\|f\|
(x-y)^2}\over {\del^2}}.$$ Since $$|f(x) - f(y)| <\eps$$ for
$|x-y|<\del$, and $$|f(x)-f(y)| < {{2\|f\| (x-y)^2}\over
{\del^2}}$$ for $|x-y|>\del$, it is readily verified that
$$p_\ell(x)\le f(x)\le p_u(x)$$ for all $x\in [a,b]$.

Since the $L_n$ are positive linear operators, this implies that
$$(L_np_\ell)(x) \le (L_nf)(x)\le (L_np_u)(x) \label{1}$$ for all
$x\in [a,b]$, and in particular for $x=y$.

For the given fixed $f$, $\eps$ and $\del$ the $p_u$ and $p_\ell$
are quadratic polynomials that depend upon $y$. Explicitly
$$p_u(x) = \left(f(y) +\eps + {{2\|f\| y^2}\over {\del^2}}\right)
- \left({{4\|f\| y}\over {\del^2}}\right) x + \left({{2\|f\|}\over
{\del^2}} \right) x^2.$$ Since the coefficients are bounded
independently of $y\in [a,b]$, and $$\lim_{n\to\infty} L_n(x^i) =
x^i,\qquad i=0,1,2,$$ uniformly on $[a,b]$, it follows that
there exists an $N$ such that for all $n\ge N$, and every choice
of $y\in [a,b]$ we have $$\left| (L_np_u)(x) -
p_u(x)\right|<\eps$$ and $$\left| (L_np_\ell)(x) -
p_\ell(x)\right|<\eps$$ for all $x\in [a,b]$. That is, $L_np_u$
and $L_np_\ell$ converge uniformly in both $x$ and $y$ to $p_u$
and $p_\ell$, respectively. Setting $x=y$ we obtain
$$(L_np_u)(y)<p_u(y) +\eps=f(y)+ 2\eps$$ and
$$(L_np_\ell)(y)>p_\ell(y) -\eps=f(y)- 2\eps.$$

Thus given $\eps>0$ there exists an $N$ such that for all $n\ge N$
and every $y\in [a,b]$ we have from (\recall{1}) $$f(y)-2\eps <
(L_nf)(y)< f(y) +2\eps.$$ This proves the theorem. \eop

A similar result holds in the periodic case $\tilC[0, 2\pi]$,
where  ``test functions'' are $1$, $\sin x$, and $\cos x$.
Numerous generalizations may be found in the book of Altomare and
Campiti [1994].

How can the Bohman-Korovkin theorem be applied to obtain density
results? It can, in theory, be applied easily. If the
$U_n=\span\{u_1\nek u_n\}$, $n=1,2,...$, are a nested sequence of
finite-dimensional subspaces of $C[a,b]$, and $L_n$ is a positive
linear operator mapping $C[a,b]$ into $U_n$ that satisfies the
conditions of the above theorem, then the $(u_k)_{k=1}^\infty$
span a dense subset of $C[a,b]$. In practice it is all too rarely
applied in this manner. The importance of the Korovkin theory is
primarily in that it presents conditions implying convergence, and
also in that it provides calculable error bounds on the rate of
approximation.

\medskip\noindent
{\bf Example 4.2.} One immediate application of the
Bohman-Korovkin theorem is a proof of the convergence of the
Bernstein polynomials $B_n(f)$ to $f$ for each $f$ in $C[0,1]$.
Recall from section 2 that for each such $f$ $$B_n(f\,;\,x) =
\sum_{m=0}^n f\left({m \over n}\right){m \choose
n}x^m(1-x)^{n-m}.$$ We can consider the $(B_n)$ as a sequence of
positive linear operators mapping $C[a,b]$ into $\Pi_n$, the space
of algebraic polynomials of degree at most $n$. It is readily
verified that $B_n(1\,;\,x) =1$, $B_n(x\,;\,x) =x$ and
$B_n(x^2\,;\,x) = x^2 + x(1-x)/n$ for all $n\ge 2$. Thus by the
Bohman-Korovkin theorem $B_n(f)$ converges uniformly to $f$ on
$[0,1]$.

\medskip\noindent
{\bf Example 4.3.} Recall from section 2 that the Fej\'er
operators $\sig_n$ maps $\tilC[0, 2\pi]$ into $T_n$. It is easily
checked that $\sig_n$ is a positive linear operator. Furthermore,
$\sig_n(1;x)=1$, $\sig_n(\sin x;x) = (n/(n+1)) \sin x$, and
$\sig_n(\cos x;x) = (n/(n+1)) \cos x$. Thus from the periodic
version of the Bohman-Korovkin theorem $\sig_n(g)$ converges
uniformly to $g$ on $[0, 2\pi]$, for each $g\in \tilC[0, 2\pi]$.

\sect{Some Univariate Density Results} {\bf Example 5.1.} {\it
M\"untz's Theorem.} Possibly the first generalization of
consequence of the Weierstrass theorems, and certainly one of the
best known, is the M\"untz theorem or the M\"untz-Sz\'asz theorem.

It was Bernstein who in a paper in the proceedings of the 1912
International Congress of Mathematicians held at Cambridge,
Bernstein [1913], and in his 1912 prize-winning essay, Bernstein
[1912], asked for exact conditions on an increasing sequence of
positive exponents $\lam_n$ so that the sequence $(x^{\lam_n})$
is fundamental in the space $C[0,1]$. Bernstein himself had
obtained some partial results. In the paper in the ICM proceedings
Bernstein wrote the following: {\it It will be interesting to know
if the condition that the series $\sum 1/\lam_n$ diverges is not
necessary and sufficient for the sequence of powers $(x^{\lam_n})$
to be fundamental; it is not certain, however, that a condition of
this nature should necessarily exist.}

It was just two years later that M\"untz [1914] was able to
provide a solution confirming Bernstein's qualified guess. What
M\"untz proved is the following.

\proclaim M\"untz's Theorem \label{}. The sequence $$x^{\lam_0},
x^{\lam_1}, x^{\lam_2}, \ldots$$ where $0\le \lam_0<\lam_1<\lam_2
<\cdots \to\infty$ is fundamental in $C[0,1]$ if and only if
$\lam_0=0$ and $$\sum_{k=1}^\infty {1\over
{\lam_k}}=\infty.\label{2}$$\nopf

There are numerous proofs and generalizations of the M\"untz
theorem. It is to be found in many of the classic texts on
approximation theory, see e.~g.~Achieser [1956, p.~43--46], Cheney
[1966, p.~193--198], Borwein, Erd\'elyi [1995, p.~171--205]. (The
last reference contains many generalizations of M\"untz's theorem
and also surveys the literature on this topic.)  We present here
the classical proof due to M\"untz, with some additions from
Sz\'asz [1916] that put M\"untz's argument into a more elegant
form.

\medskip\noindent
\pf Let $$M_n = \span\{x^{\lam_0} \nek x^{\lam_n}\},$$
and $$E(f, M_n)_\infty = \min_{p\in M_n} \|f-p\|_\infty.$$ Based
on the Weierstrass theorem it is both necessary and sufficient to
prove that $$\lim_{n\to\infty} E(x^m, M_n)_\infty =0$$ for each
$m=0,1,2,...$

To estimate $E(x^m, M_n)_\infty$ we first calculate $$E(f, M_n)_2
= \min_{p\in M_n} \|f-p\|_2,$$ where $\|\cdot\|_2$ is the
$L^2[0,1]$ norm. It is well known that $$E^2(f, M_n)_2 =
{{G(x^{\lam_0} \nek x^{\lam_n}, f)}\over {G(x^{\lam_0} \nek
x^{\lam_n})}}$$ where $G(f_1\nek f_k)$ is the Gramian of $f_1\nek
f_k$, i.e., $$G(f_1\nek f_k) = \det \left(\langle f_i\,,f_j\rangle
\right)_{i,j=1}^k.$$

As $$\langle x^p\,, x^q\rangle = \int_0^1 x^p x^q \dd x = {1\over
{p+q+1}}$$ and $$\det\left({1\over {a_i+b_j}}\right)_{i,j=1}^r
={{\prod_{1\le j<i\le r}(a_i-a_j)(b_i-b_j)}\over {\prod_{i,j=1}^r
(a_i+b_j)}},$$ a simple calculation leads to $$E^2(x^m, M_n)_2 =
{{\prod_{k=0}^n (m-\lam_k)^2}\over {(2m+1) \prod_{k=0}^n
(m+\lam_k+1)^2}}.$$ Thus, as is easily proven,
$$\lim_{n\to\infty} E(x^m, M_n)_2 =0$$ if and only if
$$\lim_{n\to\infty} \prod_{k=0}^n{{m-\lam_k}\over
{m+\lam_k+1}}=0,$$ i.e., $$\prod_{k=0}^\infty \left( 1-
{{2m+1}\over {m+\lam_k+1}}\right)=0.$$ Assuming $m\ne \lam_k$
for every $k$ (otherwise there was no reason to do this
calculation) we have $1\ne (2m+1)/(m+\lam_k+1)
>0$ and
$$\lim_{k\to\infty} {{2m+1}\over {m+\lam_k+1}}=0.$$ Thus
$$\prod_{k=0}^\infty \left( 1- {{2m+1}\over
{m+\lam_k+1}}\right)=0$$ if and only if $$\sum_{k=0}^\infty
{{2m+1}\over {m+\lam_k+1}}=\infty,$$ that in turn is equivalent
to $$\sum_{k=1}^\infty{1\over {\lam_k}}=\infty,$$ independent of
$m$. So a necessary and sufficient condition for density in the
$L^2[0,1]$ norm is that (\recall{2}) holds.

We now consider $C[0,1]$. Assume $$\sum_{k=1}^\infty{1\over
{\lam_k}}<\infty.$$ Then $E(x^m, M_n)_2$ does not tend to zero
as $n\to\infty$ for every $m$ that is not one of the $\lam_k$. As
$$E(f, M_n)_2 \le E(f, M_n)_\infty$$ for every $f\in C[0,1]$, we
have that the system $$x^{\lam_0}, x^{\lam_1}, x^{\lam_2},
\ldots$$ is not fundamental in $C[0,1]$. Furthermore, if
$\lam_0>0$ then all the functions $x^{\lam_0}, x^{\lam_1},
x^{\lam_2}, \ldots$ vanish at $x=0$, and density cannot possibly
hold.

Let us now assume that (\recall{2}) holds, and $\lam_0=0$. We will
show how to uniformly approximate each $x^m$, $m\ge 1$.  For $x\in
[0,1]$ $$\eqalign{|x^m -\sum_{k=1}^n a_k x^{\lam_k}| & = \left|
\int_0^x \left(mt^{m-1} - \sum_{k=1}^n a_k \lam_k
t^{\lam_k-1}\right)\dd t\right|\cr \le & \int_0^1 \left|mt^{m-1} -
\sum_{k=1}^n a_k \lam_k t^{\lam_k-1}\right|\dd t\cr \le &
\left(\int_0^1 \left|mt^{m-1} - \sum_{k=1}^n a_k \lam_k
t^{\lam_k-1}\right|^2 \dd t\right)^{1/2}.\cr}$$ Thus we can
approximate $x^m$ arbitrarily well in the uniform norm from the
system $x^{\lam_1}, x^{\lam_2}, \ldots$ if we can approximate
$x^{m-1}$ arbitrarily well in the $L^2[0,1]$ norm from the system
$x^{\lam_1-1}, x^{\lam_2-1}, \ldots$. We know that the latter
holds if $$\sum_{k\ge k_0} {1\over {\lam_k-1}}=\infty$$ where
$k_0$ is such that $\lam_{k_0}-1>0$. From (\recall{2}) and since
the $\lam_k$ are an increasing sequence tending to $\infty$, this
condition necessarily holds. This proves the sufficiency. \eop

The above method of showing how the $L^2$ result implies the
$C[0,1]$ result is due to Sz\'asz, and simplifies a more
complicated argument due to M\"untz that uses Fej\'er's proof of
the Weierstrass theorem. An alternative method of proof of
M\"untz's theorem and its numerous generalizations is via the
functional analytic approach, and the possible sets of uniqueness
for zeros of analytic functions, see e.~g.~Schwartz [1943], Rudin
[1966, p.~304--307], Luxemburg, Korevaar [1971], Feinerman, Newman
[1974, Chap.~X], and Luxemburg [1976]. For some different
approaches see, for example, Rogers [1981], Burckel, Saeki [1983],
and the very elegant v.~Golitschek [1983].

The above proof of M\"untz and Sz\'asz as well as most of the
functional analytic proofs, that use analytic methods, first prove
the $L^2$ result. Rudin's approach is more direct, and we
reproduce it here.

\medskip\noindent
{\bf Rudin's Proof:} Assume $0=\lam_0<\lam_1<\cdots$. If
$(x^{\lam_n})$ is not fundamental in $C[0,1]$ then from the
Hahn-Banach theorem and Riesz representation theorem there exists
a Borel measure $\mu$ of bounded total variation such that
$$\int_0^1 x^{\lam_n} \dd \mu(x)=0,$$ $n=0,1,2,\ldots$. As
$\lam_0=0$ and $\lam_n>0$ for all $n>1$ we may assume the above
holds for $n=1,2,\ldots$ and $\mu$ has no mass concentrated at
$0$. Set $$f(z) = \int_0^1 x^z \dd \mu(x).$$ For $x\in (0,1]$ and
${\rm Re}\, z>0$ we have that $x^z = e^{z \ln x}$ and $|x^z| =
x^{{\rm Re}\, z}\le 1$. It therefore follows that $f$ is analytic
and bounded in the right half plane, and of course satisfies
$$f(\lam_n)=0,\qquad n=1,2,\ldots$$ Now set $$g(z) =
f\left({{1+z}\over {1-z}}\right).$$ The transformation
$(1+z)/(1-z)$ maps the unit disc to the right half plane. Thus
$g\in H^\infty$, the space of bounded analytic functions in the
unit disc, and $g(\alp_n)=0$ where $$\alp_n = {{\lam_n-1}\over
{\lam_n+1}}.$$ Now it is a known result associated with Blaschke
products that the $(\alp_n)$ are the zeros, in the unit disc, of a
nontrivial $g\in H^\infty$ if and only if $$\sum_{n=1}^\infty
(1-|\alp_n|)<\infty.$$ It is readily checked that
$\sum_{n=1}^\infty (1-|\alp_n|)<\infty$ if and only if
$\sum_{n=1}^\infty 1/\lam_n<\infty$. Thus if $\sum_{n=1}^\infty
1/\lam_n = \infty$, then $g=0$ which implies that $f=0$. But then
$$0=f(k) = \int_0^1 x^k \dd \mu(x)$$ for all $k=1,2,\ldots$ which
implies by the Weierstrass theorem that $\mu=0$. Thus if
$\sum_{n=1}^\infty 1/\lam_n =\infty$ then the
$(x^{\lam_n})_{n=0}^\infty$ are fundamental in $C[0,1]$.

Assume that $\sum_{n=1}^\infty 1/\lam_n<\infty$. How can we
construct the desired measure $\mu$? One way is as follows. Set
$$f(z) = {1\over {(z+2)^2}}\prod_{n=0}^\infty {{\lam_n-z}\over
{2+\lam_n+z}}.$$ The function $f$ is a meromorphic function with
poles at $-2$ and $-\lam_n-2$, and zeros at the $\lam_n$. $f$ is
also bounded in ${\rm Re}\, z> -1$ since each factor is less than
1 in absolute value thereon. For each $z$ satisfying ${\rm Re}\,
z> -1$ we have by Cauchy's formula $$f(z) = {1\over {2\pi
i}}\int_{\Gamma_R} {{f(w)}\over {w-z}} dw$$ where $\Gamma_R$ is
the right semi-circle of radius $R \,(> 1+|z|)$, centered at $-1$,
together with the line from $-1-iR$ to $-1+iR$. Letting $R\to
\infty$, it may be readily shown that the integral over the
semi-circle tends to zero, and we obtain $$f(z) = {1\over {2\pi}}
\int_{-\infty}^\infty {{f(-1+is)}\over{1+z-is}}\dd s.$$ As
$${1\over{1+z-is}} = \int_0^1 x^{z-is}\dd x$$ for ${\rm Re}\,z >-1$
we have $$f(z) = \int_0^1 x^z \left({1\over {2\pi}}
\int_{-\infty}^\infty f(-1+is)e^{-is \ln x}ds\right) \dd x.$$ Set
$$\dd \mu(x) = {1\over {2\pi}} \int_{-\infty}^\infty f(-1+is)e^{-is
\ln x}\dd s.$$ This is the Fourier transform of $f(-1+is)$ at $\ln
x$ and is bounded and continuous on $(0,1]$, since the factor
$1/(2+z)^2$ in the definition of $f$ ensures that $f(-1+is)$ is a
function in $L^1$. Thus we have obtained our desired measure
$\mu$. \eop

\noindent {\bf Example 5.2.} Combining the functional analytic
approach with analytic methods has proven to be a very effective
method of proving density results. As a general example, assume
$g$ is in $C(\RR)$ and has an extension as an analytic function on
all of $\CC$. Let $\Lam$ be a subset of $\RR$ that contains a
finite accumulation point, i.e., there are distinct $\lam_n$ in
$\Lam$ and a finite $\lam^*$ such that $\lim_{n\to\infty}
\lam_n=\lam^*$. Set $$\calM_\Lam =\span \{ g(\lam x):\, \lam\in
\Lam\}.$$ We wish to determine when $\calM_\Lam$ is dense in
$C[a,b]$. The following result holds.

\proclaim Theorem \label{5.2}. Let $g$, $\Lam$ and $\calM_\Lam$ be
as above. Set $$N_g =\{ n:\, g^{(n)}(0)\ne 0 \}.$$ Then
$\calM_\Lam$ is dense in $C[a,b]$ if and only if:
\item{i)} for $[a,b] \incl (0,\infty)$ or $[a,b] \incl (-\infty, 0)$
$$\sum_{n\in N_g\\ \{0\}}{1\over n} =\infty,$$
\item{ii)} if $a=0$ or $b=0$, then $0\in N_g$ and
$$\sum_{n\in N_g\\ \{0\}}{1\over n} =\infty,$$
\item{iii)} if $a<0<b$, then $0\in N_g$ and
$$\sum_{{{n\in N_g\\ \{0\}} \atop {n\ {\rm even}}}} {1\over n} =
\sum_{{{n\in N_g} \atop {n\ {\rm odd}}}} {1\over n} =\infty.$$

\pf  The conditions in (i), (ii) and (iii) are exactly
those conditions that determine when $$\span\{x^n :\, n\in N_g\}$$
is dense in $C[a,b]$. This is the content of the M\"untz theorem
in case (ii), and easily follows from the M\"untz theorem in case
(iii). In case (i) it follows from the M\"untz theorem that the
condition therein is sufficient for density. The
necessity is also true, but needs an additional argument, see
e.g., Schwartz [1943].

From the Hahn-Banach and Riesz representation theorems
$\calM_\Lam$ is not dense in $C[a,b]$ if and only if there exists
a nontrivial measure $\mu$ of bounded total variation on $[a,b]$
satisfying $$\int_a^b g(\lam x)\dd \mu(x) =0$$ for all $\lam\in
\Lam$. Assume such a measure exists. As $g$ is entire, it follows
that $$h(z) = \int_a^b g(z x)\dd \mu(x)$$ is entire. Furthermore
$h(\lam)=0$ for all $\lam\in \Lam$. By assumption $\Lam$ contains
a finite accumulation point. Thus by the uniqueness theorem for
zeros of analytic functions $h=0$. However $h$ being identically
zero does not necessarily imply that $\mu$ is the zero measure. It
only proves that $$\overline{\calM_\Lam} =\overline{\span} \{
g(\lam x):\, \lam\in \RR\}.$$ For example, if $g$ is a
polynomial of degree $m$, then $\calM_\Lam$ is simply the space of
polynomials of degree $m$.

As $$\int_a^b g(z x)\dd \mu(x) =0$$ and $g$ is entire it may be
shown, differentiating by $z$, that $$\int_a^b x^n g^{(n)}(z
x)\dd \mu(x) =0$$ for every nonnegative integer $n$. Setting $z=0$
gives us $$g^{(n)}(0)\int_a^b x^n \dd \mu(x) =0,\qquad n=0,1,...$$
Thus $$\int_a^b x^n \dd \mu(x) =0,$$ for all $n\in N_g$. But
$\span\{ x^n:\, n\in N_g\}$ is dense in $C[a,b]$, so $\mu$ is the
trivial measure.

On the other hand, assume the conditions in (i), (ii) or (iii) do
not hold. Thus $\span\{ x^n:\, n\in N_g\}$ is not dense in
$C[a,b]$, and there exists a nontrivial measure $\mu$ of bounded
total variation satisfying $$\int_a^b x^n \dd \mu(x) =0$$ for all
$n\in N_g$. Since $g$ is entire $$g(x) = \sum_{n\in N_g}
{{g^{(n)}(0)}\over {n!}} x^n$$ and it follows that $$\int_a^b
g(\lam x) \dd \mu(x) =0$$ for all $\lam\in \RR$. $\calM_\Lam$
is not dense in $C[a,b]$. \eop

For example, if $g(x)=e^x$ then $N_g = \ZZ_+$ so that (i), (ii)
and (iii) always hold. Thus $$\span\{ e^{\lam_n x}:\, \lam\in \Lam
\}$$ is always dense in $C[a,b]$ assuming $\Lam$ is a subset of
$\RR$ with a finite accumulation point. A change of variable
argument implies that under this same condition on $\Lam$ the set
$$\span\{ x^{\lam_n}:\,\lam \in \Lam\}$$ is dense in $C[\alp,
\beta]$ for every $0<\alp<\beta<\infty$.

A question related to M\"untz type problems is that of the
fundamentality of the functions $(e^{\lam_n x})$, where $(\lam_n)$
is a sequence of complex numbers. This has been considered in the
space of complex-valued functions in $C[a,b]$, $C(\RR_+)$,
$L^p[a,b]$ and $L^p(\RR_+)$, $1\le p<\infty$. There has been a
great deal of research done in this area, see, for example, Paley,
Wiener [1934, Chap.~VI], Levinson [1940, Chap.~I and II], Schwartz
[1943], Levin [1964, Appendix III], Levin [1996, Lecture 18], and
the many references therein.

\medskip\noindent
{\bf Example 5.3.} {\it Akhiezer's Theorem.} Let $\Gamma$ be a
subset of $\RR \\ [-1,1]$, and consider the set $$\calN_\Gamma
=\span \left\{{1\over {t-\gam}}:\,\gam\in \Gamma\right\}.$$ When
is $\calN_\Gamma$ dense in $C[-1,1]$? One result is similar to
Theorem \recall{5.2}. It may be found in  Feinerman-Newman [1974,
p.~116--117], but the proof therein is somewhat different.

\proclaim Proposition \label{5.3}. If $\Gamma$ has either a finite
accumulation point in $\RR \\ [-1,1]$ or $\infty$ is an
accumulation point, then $\calN_\Gamma$ is dense in $C[-1,1]$.

\pf Assume $\overline{\calN_\Gamma}\ne C[-1,1]$. Then
there exists a Borel measure $\mu$ of bounded total variation such
that $$\int_{-1}^1  {1\over {t-\gam}} \dd \mu(t) =0$$ for all
$\gam\in \Gamma$. Set $$f(z) = \int_{-1}^1  {1\over {t- z}}
\dd \mu(t).$$ Note that $f(\gam)=$ for all $\gam\in \Gamma$. It
is readily verified that $f$ is analytic on $\CC \\ [-1,1]$, and
analytic also at infinity.

Thus if  $\Gamma$ has either a finite accumulation point in $\RR
\\ [-1,1]$ or $\infty$ is an accumulation point, then $f=0$. For $|z|>1\ge
|t|$ $$f(z) = -\sum_{n=0}^\infty \left({1\over z}\right)^{n+1}
\int_{-1}^1 t^n \dd \mu(t).$$ As $f=0$, this then implies that
$$\int_{-1}^1 t^n \dd \mu(t)=0$$ for all $n$, which from the
Hahn-Banach and Weierstrass theorems implies that $\mu$ is the
trivial measure. This proves the proposition. \eop

What can be said if the only accumulation points of $\Gamma$ are
$1$ or $-1$ or both? The result is known, contains the previous
Proposition \recall{5.3} as a special case, and was proven by
Akhiezer, see Achieser [1956, p.~254--256].

\proclaim Akhiezer's Theorem \label{5.4}. Let
$(\gam_n)_{n=1}^\infty$ be a sequence in $\RR \\ [-1,1]$, and
consider the set $$\calN =\span \left\{{1\over {t-\gam_n}}:\,
n=1,2,\ldots\right\}.$$ Then $\calN$ is dense in $C[-1,1]$ if
and only if $$\sum_{n=1}^\infty 1 - |\gam_n - \sqrt{\gam_n^2
-1}|=\infty.$$\nopf

See also Borwein, Erd\'elyi [1995, p.~208] where a different
method of proof is used. They also give the above condition as
$$\sum_{n=1}^\infty \sqrt{\gam_n^2 -1}=\infty,$$ and these two
conditions are in fact equivalent. Akhiezer's proof of this theorem
is delicate and detailed, dependent on the construction of
specific best approximants. We will not reproduce it here. Michael
Sodin has a proof which uses complex variable theory.

\medskip\noindent
{\bf Example 5.4.} The analysis literature is replete with results
concerning the density of {\it translates (and dilates)} of a
function in various spaces. These might be arbitrary, integer, or
sequence translates (or dilates). Many of these results are
generalizations, in a sense, of the M\"untz and/or Paley-Wiener
theorems. See, for example, both Example 5.2 and 5.3.

There is a characterization of those $f\in C(\RR)$ for which
$$\span\{ f(\cdot\,-\alp)\,:\, \alp\in \RR\}$$ is {\bf not} dense
in $C(\RR)$ (in the topology of uniform convergence on compacta).
Such functions are called {\it mean periodic}, see Schwartz
[1947].

Some functions in $C(\RR)$ have a further interesting property.

\proclaim Proposition \label{5.5}. Assume $f=\widehat{g}$ ($f$ is
the Fourier transform of $g$) for some nontrivial $g\in L^1(\RR)$
with the support of $g$ contained in an interval of length at most
$2\pi$. Then $$\span\{ f(\cdot\,-n)\,:\, n\in \ZZ\}$$ is dense in
$C(\RR)$ (in the topology of uniform convergence on compacta).

\pf Assume the above set is not dense in $C(\RR)$. There
then exists a Borel measure $\mu$ of bounded total variation and
compact support $E$ such that $$\int_E  f(x-n) \dd \mu(x) =0$$ for
all $n\in \ZZ$. Assume $f=\widehat{g}$, as above, and
$\supp\{g\}\incl [a, a+2\pi]$. Thus for each $n\in \ZZ$ $$0=
\int_E  f(x-n) \dd \mu(x) = \int_E  \widehat{g}(x-n) \dd \mu(x)$$
$$= {1\over {2\pi}} \int_E\left( \int_a^{a+2\pi}
g(t)e^{-i(x-n)t}\dd t\right) \dd \mu(x)$$ $$= \int_a^{a+2\pi}
\left({1\over {2\pi}}\int_E e^{-ixt} \dd \mu(x)\right)
e^{int}g(t)\dd t =\int_a^{a+2\pi}
e^{int}g(t)\widehat{\mu}(t)\dd t$$
where $\widehat{\mu}$ is the
Fourier transform of the measure $\mu$. It is well known that
$\widehat{\mu}$ is an entire function.

As all the Fourier coefficients of $g\, \widehat{\mu}$ on $[a,
a+2\pi]$ vanish we have that $g\,\widehat{\mu}$ is identically
zero thereon. This implies that $g$ must vanish where
$\widehat{\mu}\ne 0$. As $\widehat{\mu}$ is entire this implies
that $g=0$, a contradiction. \eop

The above is a simple example within a general theory. The
interested reader should consult Atzmon, Olevskii [1996], Nikolski
[1999], and references therein. Note that there is no function
whose integer translates are dense in $L^2(\RR)$.

\medskip\noindent
{\bf Example 5.5.} {\it The Bernstein Approximation Problem.}
Assume $\ome$ is a {\it weight} on $\RR$ by which we will mean a
non-negative, measurable, bounded function. For each $f$ in
$C(\RR)$ satisfying $$\lim_{|x|\to \infty} \ome(x)f(x) =0$$ set
$$\|f\|_\ome = \sup_{x\in \RR} \ome(x)|f(x)|,$$ and let
$C_\ome(\RR)$ denote the real normed linear space of those $f$ as
above with $\|f\|_\ome <\infty$. The Bernstein approximation
problem was first formulated in Bernstein [1924]. It asks for
necessary and sufficient conditions on a weight $\ome$ such that
(algebraic) polynomials are dense in $C_\ome(\RR)$. That is, for
each $f$ in $C_\ome(\RR)$ and $\eps>0$ there exists a polynomial
$p$ for which $\|f-p\|_\ome <\eps$. This immediately implies that
$\ome$ must satisfy $$\lim_{|x|\to\infty} \ome(x)p(x) =0$$ for
every polynomial $p$.

In Bernstein [1924] can be found the following result. Assume
$\ome(x)=1/q(x)$, where $$q(x) =\sum_{n=0}^\infty a_n x^{2n}$$
with $a_0>0$, $a_n \ge 0$ for all $n$, and $q$ not the constant
function. Then a necessary and sufficient for polynomials to be
dense in $C_\ome(\RR)$ is that $$\int_1^\infty {{\ln q(x)}\over {1
+ x^2}}\dd x =\infty.$$ In general a condition of this form is
necessary, but not sufficient. It is often sufficient for
``reasonable'' weights.

The literature on this problem is rather extensive including the
review articles Ahiezer [1956] and Mergelyan [1956], see also
Lorentz, v.~Golitschek and Makovoz [1996, p.~28-33], and Timan
[1963, p.~16--19]. The article of Mergelyan, as well as Prolla
[1977], includes a proof of this next result. Let $\calM_\ome$
denote the set of polynomials $p$ satisfying $\ome(x)|p(x)| \le 1
+ |x|$ for all $x\in \RR$, and set $$M_{\ome}(z) =\sup\{|p(z)|:\,
p\in \calM_\ome\}.$$

\proclaim Mergelyan's Theorem \label{5.6}. Let $\ome$ be as above.
Then a necessary and sufficient for polynomials to be dense in
$C_\ome(\RR)$ is that $$M_{\ome}(z) =\infty$$ for every $z\in
\CC\\ \RR$.\nopf

Unfortunately this condition is not easy to check.

Here is a condition that is easier to check, but which only holds
for certain weights. Assume $\ome =\exp \{-Q\}$, $\ome$ is even,
and $Q$ is a convex function of $\ln x$ on $(0,\infty)$. Then
$$\int_0^\infty {{\ln \ome(x)}\over {1 + x^2}}\dd x=-\infty$$ is both
necessary and sufficient for polynomials to be dense in
$C_\ome(\RR)$, see Mhaskar [1996, p.~331].

\medskip\noindent
{\bf Example 5.6.} {\it Markov Systems.} Assume we are given a
sequence of functions $(u_m)_{m=0}^\infty$ in $C[a,b]$. What we
have been asking is when this sequence is fundamental, i.e.,
linear combinations are dense. That is when, for each $f$ in
$C[a,b]$ and $\eps>0$,  there exists a finite linear
combination $u$ of the $(u_m)_{m=0}^\infty$ such that
$$\|f-u\|_\infty < \eps.$$ The only general result
characterizing the density of such sequences is the somewhat
tautological Theorem \recall{banach2}. However, when the sequence
$(u_m)_{m=0}^\infty$ has a particular type of Chebyshev property,
then P.~Borwein proved a surprisingly interesting condition
equivalent to density.

To explain his result we first need some definitions. Given the
sequence $(u_m)_{m=0}^\infty$ in $C[a,b]$ we set $$U_n =\span\{
u_0,u_1\nek u_n\}$$ for each $n\in \ZZ_+$. We say that $U_n$ is a
{\it Chebyshev space} if no $u\in U_n$, $u\ne 0$, has more than
$n$ distinct zeros in $[a,b]$. We say that the sequence
$(u_m)_{m=0}^\infty$ is a {\it Markov sequence} if $U_n$ is a
Chebyshev space for $n=0,1\ldots$. There are numerous examples of
Markov sequences. For example, $(x^{\lam_m})_{m=0}^\infty$ is a
Markov sequence on $[a,b]$ where $a>0$ and the
$(\lam_m)_{m=0}^\infty$ are arbitrary distinct real values, while
$(1/(x-c_m))_{m=0}^\infty$ is a Markov sequence on any $[a,b]$
where the $(c_m)_{m=0}^\infty$ are distinct values in $\RR\\
[a,b]$.

In what follows we assume that the $(u_m)_{m=0}^\infty$ is a
Markov sequence. Let $t_n\in U_n$ be of the form $t_n = u_n - v_n$
with $v_n\in U_{n-1}$, satisfying $$\|t_n\|_\infty = \min_{v\in
U_{n-1}}\|u_n - v\|_\infty.$$ It is well known, from the
Chebyshev and Markov properties, that $t_n$ is uniquely defined
and has $n$ zeros in $(a,b)$. Let $x_1<\cdots < x_n$ denote these
$n$ zeros and set $x_0=a$, $x_{n+1}=b$. The {\it mesh} of $t_n$ is
defined by $$m_n =\max_{i=1\nek n+1} (x_i - x_{i-1}).$$ It is
readily proven that for any $k<n$ the function $t_k$ has at most
one zero between any two consecutive zeros of $t_n$. From this it
follows that $$\lim_{n\to\infty} m_n =0$$ if and only if
$$\liminf_{n\to\infty} m_n=0.$$

The following result may be found in Borwein [1990], and also in
Borwein, Erd\'elyi [1995, p.~155--158].

\proclaim Borwein's Theorem \label{5.7}. Assume
$(u_m)_{m=0}^\infty$ is a Markov sequence in $C^1[a,b]$ and
$u_0=1$. Then the sequence $(u_m)_{m=0}^\infty$ is dense in
$C[a,b]$ if and only if
$$\lim_{n\to\infty} m_n =0.$$\nopf

A similar result relating density to Bernstein-type inequalities
is in Borwein, Erd\'elyi [1995b], and Borwein, Erd\'elyi [1995,
p.~206--211].

\medskip\noindent
{\bf Example 5.7.} The following result is a special case of a
general theorem of Schwartz [1944] (see also Pinkus [1996] and
references therein). Here we again consider $C(\RR)$, with the
topology of uniform convergence on compacta. We are interested in
determining the set of functions in $C(\RR)$ that are both
translation and dilation invariant.

\proclaim Proposition \label{5.8}. If $\sig\in C(\RR)$, $\sig\ne
0$, then
$$C(\RR)=\overline{\span}\{\sig(\alp \cdot +\beta):\alp, \beta\in
\RR\}$$ if and only if $\sig$ is not a polynomial.

\pf Let $$\calM_\sig ={\span}\{\sig(\alp
\cdot+\beta):\alp, \beta \in \RR\}.$$ If
$\overline{\calM_\sig}\ne C(\RR)$ then there exists a nontrivial
Borel measure $\mu$ of bounded total variation and compact support
such that $$\int_{\RR} \sig(\alp x+\beta)\dd \mu(x)=0$$ for all
$\alp, \beta\in \RR$. Since $\mu$ is nontrivial and polynomials
are dense in $C(\RR)$ in the topology of uniform convergence on
compact subsets, there must exist a $k \ge 0$ such that
$$\int_{\RR} x^k\dd \mu(x)\ne 0.$$

It is relatively simple to show that for each $\phi\in
C^\infty_0(\RR)$, (infinitely differentiable and having compact
support) the convolution $(\sig*\phi)$ is contained in
$\overline{\calM_\sig}$. Since both $\sig$ and $\phi$ are in
$C(\RR)$, and $\phi$ has compact support, this can be proven by
taking limits of Riemann sums of the convolution integral. We also
consider taking derivatives as a limiting operation in taking
divided differences. Since $(\sig*\phi)\in C^\infty(\RR)$, and
thus it and all its derivatives are uniformly continuous on every
compact set, it follows that for each $\alp, \beta\in \RR$
$${{\partial^n}\over {\partial \alp^n}} (\sig*\phi)(\alp x+\beta)
= x^n (\sig*\phi)^{(n)}(\alp x+\beta)\in \overline{\calM_g}.$$
Thus $$\int_{\RR} x^n (\sig*\phi)^{(n)}(\alp x+\beta)
\dd \mu(x)=0,$$ for all $\alp, \beta \in \RR$ and $n\in \ZZ_+$.
Setting $\alp=0$, we see that $$(\sig*\phi)^{(n)}(\beta)\int_{\RR}
x^n \dd \mu(x)=0$$ for each choice of $\beta\in \RR$, $n\in \ZZ_+$
and $\phi\in C_0^\infty(\RR)$. This implies, since $\int_{\RR}
x^k\dd \mu(x)\ne 0$, that $$(\sig*\phi)^{(k)}=0$$ for all $\phi\in
C_0^\infty(\RR)$. That is, $\sig^{(k)}=0$ in the weak sense.
However, as is well-known, this implies that $\sig^{(k)}=0$ in the
strong (usual) sense. That is, $\sig$ is a polynomial of degree at
most $k-1$.

The converse direction is simple. If $\sig$ is a polynomial of
degree $m$, then $\calM_\sig$ is exactly the space of polynomials
of degree $m$, and is therefore not dense in $C(\RR)$.\eop

\noindent {\bf Example 5.8.} {\it Splines} are piecewise
polynomials with a high order of continuity. For $a= \xi_0
<\xi_1<\cdots< \xi_k < \xi_{k+1}=b$ we set
$$\calS_{n}(\xi_1\nek \xi_k) = \{ s\in C^{(n-1)}[a,b]:\, s\in \Pi_n|_{(\xi_{i-1},
\xi_i)},i=1\nek k+1\}.$$
Hence a function $s$ belongs to
$\calS_{n} (\xi_1\nek \xi_k)$ if it is a $C^{(n-1)}$ function, i.e.,
has a certain global level of
smoothness, and is a polynomial of degree at most $n$ on each of
the intervals $(\xi_{i-1}, \xi_i)$. We say that $\calS_{n}
(\xi_1\nek \xi_k)$ is the space of {\it splines of degree $n$ with
the simple knots} $(\xi_1\nek \xi_k)$. When using splines one
fixes the degree and permits the number (and placement) of the
knots to vary. From the perspective of numerical computations,
approximation by splines enjoys many advantages over approximation
by algebraic and trigonometric polynomials. As the number of knots
increases the corresponding space of splines may or may not
``become dense'' in $C[a,b]$. Whether it does or not simply
depends upon if the knots become dense in $[a,b]$.

To be more exact, for each $k=1, 2,...$ let $$S_k
=\calS_{n}(\xi_1^k\nek \xi_k^k)$$ for some set of $k$ knots as
above, where $\xi^k_0=a$ and $\xi^k_{k+1}=b$. For each such $k$,
let $$m_k = \max_{i=0\nek k} (\xi^k_{i+1} - \xi^k_i),$$ denote
the maximum mesh length. Then we have, see for example, de Boor
[1968],

\proclaim Proposition \label{5.9}. For each $f\in C[a,b]$ there
exist $s_k\in S_k$ such that $$\lim_{k\to\infty}
\|f-s_k\|_\infty=0$$ if and only if $\lim_{k\to\infty} m_k=0$.

\pf We first assume that  $\lim_{k\to\infty} m_k=0$. The
set $$(1,x\nek x^n, (x-\xi_1^k)_+^n\nek (x-\xi_k^k)_+^n),$$
where $x_+^n$ equals $x^n$ for $x\ge 0$ and $0$ for $x < 0$, is a
basis for $S_k$. However there are also  `better' bases. They are
given by B-splines $$(B_1^k\nek B_{n+k+1}^k).$$ Each $B_i^k$ is
nonnegative, $\sum_{i=1}^{n+k+1}B_i^k =1$ on $[a,b]$, and
$\overline{\supp}\, B_i^k = [\xi_{i-n-1}^k, \xi_i^k]$ where
$\xi_{-n}^k \le \cdots \le \xi_{-1}^k\le a$, and $b\le
\xi_{k+2}^k\le \cdots \le \xi_{n+k+1}^k$.

For each $\del>0$, let $$\ome(f;\del) = \max_{|x-y|\le \del}
|f(x)-f(y)|$$ denote the usual {\it modulus of continuity} of $f$.
As $f$ is continuous on $[a,b]$ it is also uniformly continuous
thereon and thus $$\lim_{\del \to 0^+} \ome(f;\del)=0.$$

Now choose $t_i^k\in \supp\, B_i^k\cap [a,b]$, $i=1\nek n+k+1$,
and set $$s_k(x) =\sum_{i=1}^{n+k+1} f(t_i^k) B_i^k(x).$$ ($s_k$
is called a {\it quasi-interpolant}.) Then for $x\in [\xi_j^k,
\xi_{j+1}^k]$, $j\in \{0,1\nek k\}$ we have $$\eqalign{
|f(x)-s_k(x)| & = |f(x) - \sum_{i=1}^{n+k+1} f(t_i^k) B_i^k(x)|\cr
\le & \sum_{i=1}^{n+k+1} |f(x) - f(t_i^k)| B_i^k(x)\cr = &
\sum_{i=j}^{n+j+2} |f(x) - f(t_i^k)| B_i^k(x)\cr}$$ since
$B_i^k(x) =0 $ for $i\notin\{j\nek n+j+2\}$. Thus $$|f(x)-s_k(x)|
\le \max_{i=j\nek n+j+2} |f(x) - f(t_i)| \le \ome(f;
(n+2)m_k),$$ and $$\|f-s_k\| \le \ome(f; (n+2)m_k),$$ from
which we obtain $$\lim_{k\to\infty} \|f-s_k\|=0.$$

If the $m_k$ do not tend to zero, then there is a subinterval
$[c,d]$ of $[a,b]$ of positive length, and a subsequence $(k_m)$
of $(k)$ such that $S_{k_m}$ has no knots in $[c,d]$. That is,
each function in $S_{k_m}$ is a polynomial of degree at most $n$
on $[c,d]$. If $f\in C[a,b]$ is not a polynomial of degree at most
$n$ on $[c,d]$, then there exists a $C>0$ such that $$\|f-s\|\ge
C$$ for all $s\in S_{k_m}$.\eop

Each $\calS_n(\xi_1\nek \xi_k)$ is a linear space of splines of
degree $n$ with $k$ fixed knots. One can also consider the
nonlinear set $\calS_{n,k}$ of splines of degree $n$ with $k$ {\it
free knots}, i.e., $$\calS_{n,k} = \{ s(x) =\sum_{i=0}^n c_ix^i
+\sum_{i=1}^k d_i (x-\xi)_+^n\,:\, c_i, d_i \in \RR,\, a<
\xi_1<\cdots <\xi_k<b\}.$$ (Note that the set $\calS_{n,k}$ is
not closed.)

The following is a consequence of our previous density result.

\proclaim Proposition \label{5.10}. Let $p\in (1,\infty)$. Assume
$f\in L^p[a,b] \\ \overline{\calS_{n,k}}$ and $s^*\in
\overline{\calS_{n,k}}$ satisfies $$\|f-s^*\|_p \le \|f-s\|_p$$
for all $s\in \calS_{n,k}$. Then $s^*\notin
\overline{\calS_{n,k-1}}$.

\pf Let $F$  be a
continuous linear functional on $L^p[a,b]$ that satisfies
$$\|F\|=1$$ and $$F(f-s^*)=\|f-s^*\|\,.$$ As $f-s^*\ne 0$ and
$p\in (1,\infty)$, such an $F$ exists and is unique.

Our proof is by contradiction. If $s^*\in \overline{\calS_{n,k-1}}$,
then for each $d\in \RR$ and
$\xi\in [a,b]$ we have $$\|f-s^*\|_p \le \|f-s^*
-d(x-\xi)_+^n\|_p.$$ Thus from Theorem \recall{banach1} we have
$$F((x-\xi)_+^n)=0$$ for each $\xi\in [a,b]$. In addition, as
$$\|f-s^*\|_p \le \|f-s^* -p\|_p$$ for every $p\in \Pi_n$, we have
$$F(x^k)=0,\qquad k=0\nek n.$$ From Proposition \recall{5.9} $$\span
\{1, x\nek x^n, (x-\xi)_+^n\,:\, \xi\in [a,b]\}$$ is dense in
$C[a,b]$ and thus in $L^p[a,b]$. From Theorem \recall{banach2}
this implies that $F=0$. A contradiction. \eop

The exact same argument proves, for example, that if $$\calP_k =\{
\sum_{j=1}^k a_j x^{m_j}:\, a_j\in \RR, m_j\in \ZZ_+, j=1\nek
k\},$$ then in $L^p[a,b]$, $1<p<\infty$, any best approximation
to $f$ from $\calP_k$ is never contained in $\calP_{k-1}$.

The key ingredients of the above argument are the density of the
set under consideration (in the above examples, splines and
algebraic polynomials) in the normed linear space $E$, and the
fact that $E$ is {\it smooth}. That is, to each nonzero element of
the space $E$ there is a unique continuous linear functional of
norm one that attains its norm on the given element.

\medskip\noindent
{\bf Example 5.9.} Here are two examples where we consider the
density of positive cones. That is, we present some applications
of Theorem \recall{3.8}.

Let $\Pi$ denote the space of all algebraic polynomials and
$\Pi_+$ the positive cone of all algebraic polynomials with
nonnegative coefficients. We first prove the following result due
to Bonsall [1958].

\proclaim Theorem \label{5.11}. The uniform closure of $\Pi_+$ on
$[-1,0]$ is exactly the set of $f$ in $C[-1,0]$ for which $f(0)\ge
0$.

\pf  Let $g_n(x)= (1+x)^n$ and $\phi_n(x)=x^n$ for all
$n\in \ZZ_+$. Note that $g_n$ and $\phi_n$ are in $\Pi_+$, and
$$g_n = \sum_{k=0}^n {n \choose k} \phi_k.$$ Assume $F$ is a
continuous linear functional on $C[-1,0]$ satisfying
$F(\phi_n)\ge 0$ for every $n\in \ZZ_+$. Then $$F(g_n) \ge {n
\choose k} F(\phi_k).$$ Now $\|g_n\|=\|\phi_n\|=1$ for all $n$.
Thus $$\|F\|\ge {n \choose k} F(\phi_k),$$ for each $n\ge k$.
Fix $k\ge 1$ and let $n\to \infty$. This implies that
$F(\phi_k)=0$ for all $k=1,2,...$. Thus each $\pm \phi_k$, $k\ge
1$, is in the uniform closure of $\Pi_+$ on $[-1,0]$. As every
$f\in C[-1,0]$ satisfying $f(0)=0$ is in the uniform closure of
the space generated by the $\pm \phi_k$, $k\ge 1$, the result now
easily follows.

Bonsall actually proves that each $F$ as above is necessarily of
the exact form $F(f) = c f(0)$ where $c=F(1)\ge 0$. This he proves
as follows. For each $f\in C[-1,0]$ and $\eps>0$, let $p\in \Pi$
satisfy $$\|f-p\|<\eps$$ Thus we have $|f(0) -p(0) |<\eps$. Since
$F(\phi_k) = 0$, $k\ge 1$, we obtain $$ |F(f) - F(1)p(0)| =
|F(f-p(0))|= |F(f-p)| < \eps \|F\|.$$ Furthermore $$|F(f) -
F(1)f(0)| = |F(f) -F(1)p(0) + F(1)p(0) - F(1)f(0)| < 2\eps
\|F\|.$$ As this is valid for each $\eps>0$ we obtain $$F(f) =
F(1)f(0).\meop$$

There is an alternative method of proving this result via a slight
generalization of the Stone-Weierstrass theorem.
Consider the set of $f$ in $C[-1,0]$ satisfying
$f(0)=0$. Now $e^{\alp x}-1$ is in the uniform closure of $\Pi_+$
on $[-1,0]$ for $\alpha>0$. (Truncate the power series expansion
about 0.) Furthermore if $f$ and $g$ are in this closure then so
is $fg$. As $e^{\alp x}-1$ approaches $-1$ uniformly on $[-1,
\delta]$ for any $\delta<0$ as $\alpha\to\infty$, and is bounded
on $[\delta, 0]$, it follows that for $f$ in the uniform closure
of $\Pi_+$ on $[-1,0]$ and satisfying $f(0)=0$
we also have $-f$ in this same closure. In
addition $p(x)=x$ is nonzero for all $x\ne 0$ and separates
points. Thus from an elementary generalization of the
Stone-Weierstrass theorem the uniform closure of $\Pi_+$ on
$[-1,0]$ contains the set of all $f$ in $C[-1,0]$ satisfying $f(0)=0$.
The result now follows.

\medskip
Thus for any $a<b<0$ the uniform closure of $\Pi_+$ on $[a,b]$ is
exactly all of $C[a,b]$. What happens if $[a,b]\incl [0,\infty)$?
It is well known that in this case the uniform closure of $\Pi_+$
is simply the set of analytic functions in $[a,b]$ given by a
power series about 0 with nonnegative
coefficients which converges in $[a,b]$.

There are also somewhat surprising results due to Nussbaum, Walsh
[1998], generalizing work of Toland [1996]. These results are used
to investigate when the spectral radius of a positive, bounded
linear operator belong to its spectrum. A special case of what
they prove is the following:

\proclaim Theorem \label{5.12}. For any $a<-1$ the uniform closure
of $\Pi_+$ on $[a,1]$ contains the set of all $f$ in $C[a,1]$ that
vanishes identically on $[-1,1]$.

\pf We present two proofs of this result. The first proof
uses the Hahn-Banach theorem and is that found in Nussbaum, Walsh
[1998]. The second proof is constructive.

Assume we are given any continuous linear functional $F$ on
$C[a,1]$ satisfying $F(x^n)\ge 0$ for all $n\in \ZZ_+$. From the
Riesz representation theorem, this implies the existence of a
Borel measure $\mu$ of bounded total variation satisfying
$$\int_a^1 x^n \dd \mu(x) \ge 0$$ for all $n\in \ZZ_+$. We will prove
that $\supp \{\mu\}\incl [-1,1]$. As this is true then $$\int_a^1
f(x) \dd \mu(x)=0$$ for every $f$ in $C[a,1]$ that vanishes
identically on $[-1,1]$, proving our theorem.

To this end, consider $$G(z) = \int_a^1 {1\over{ z-x}}\dd \mu(x).$$
$G$ is analytic in $\CC \\ [a,1]$, and vanishes at $\infty$. For
$|z| >\lam = \sup\{|x|:\, x\in \supp \{\mu\}\}$ we have $$G(z)
=\sum_{n=1}^\infty {{c_{n-1}}\over {z^n}}$$ where $$c_n = \int_a^1
x^n \dd \mu(x) \ge 0.$$ Note that $H(z) = G(1/z)$ is analytic in
$\CC\\ \{(-\infty, 1/a]\union [1,\infty)\}$ and has about the
origin a power series expansion with nonnegative coefficients.
From a theorem of Pringsheim, if the radius of convergence of the
power series is $\rho>0$ then the point $z=\rho$ is a singular
point of the analytic function represented by the power series. As
the power series converges on $[0,1)$ the radius of convergence is
at least 1, and therefore $H$ is analytic in  $\CC\\
\{(-\infty,-1]\union [1,\infty)\}$ and $G$ analytic in $\CC \\
[-1,1]$. That is, $G$ is in fact analytic in $[a, -1)$. This
implies, see Nussbaum, Walsh [1998, p.~2371], that the measure
$\dd \mu$ has no support in $[a, -1)$.

The following constructive proof of this result is based on a
variation of a proof to be found in Orlicz [1992, p.~99]. For
$n\in \NN$, odd, consider the function $$g_n(x) =\int_0^x
e^{t^n/n} - 1 \dd t.$$ Note that the integrand is uniformly
bounded on $[a,1]$ and $$\lim_{n\to \infty} e^{t^n/n} - 1 =
\cases{0, & $-1\le t \le 1$\cr -1, & $a\le t < -1\,.$\cr}$$ As
$$e^{t^n/n} -1 =\sum_{k=1}^\infty \left({{t^n}\over n}\right)^k
{1\over {k!}}$$ this function is in the uniform closure of
$\Pi_+$. Thus, so is $g_n$ . Set $$G(x) =\cases{-(x+1), & $a\le
x\le -1$\cr 0, & $-1\le x \le 1\,.$\cr}$$ Then $$\lim_{n\to
\infty} (G(x) - g_n(x)) =0$$ uniformly on $[a,1]$. That is, for
all $x\in [-1,1]$ $$\left|\int_0^x  e^{t^n/n} - 1 \dd t\right| \le
e^{1/n} -1,$$ while for $x\in [a, -1]$
$$\left| -(x+1) -
\int_0^x  \left(e^{t^n/n} - 1\right) \dd t \right|$$ $$ = \left| -(x+1)
- \int_{-1}^x  \left(e^{t^n/n} - 1\right) \dd t + \int_{-1}^0
\left(e^{t^n/n} - 1\right) \dd t \right|$$ $$ \le \int^{-1}_x
e^{t^n/n}\dd t + \int_{-1}^0  \left(1- e^{t^n/n}\right) \dd t \le
\int^{-1}_a  e^{t^n/n}\dd t + (1 - e^{-1/n}).$$ Thus $G$ is in the
uniform closure of $\Pi_+$ on $[a,1]$.

Moreover, as seen above, the function $ e^{t^n/n} - 1$ is
uniformly bounded and approaches $$H(x) = \cases{0, & $-1\le x\le
1$ \cr -1, & $a\le x <-1.$\cr}$$ The convergence to $H$ is uniform
in $[a,1]$, away from any neighbourhood of $-1$. Thus $GH = -G$ is
also in the uniform closure of $\Pi_+$ on $[a,1]$, and therefore
the uniform closure of $\Pi_+$ on $[a,1]$ contains the algebra
generated by $G$. An elementary generalization of the
Stone-Weierstrass theorem implies that the uniform closure of
$\Pi_+$ on $[a,1]$ contains the set of all $f$ in $C[a,1]$  which
vanish identically on $[-1,1]$. \eop

The above result is an extension of Theorem II$'$ in Orlicz [1999,
p.~96]. Orlicz proved that for every $f\in C[a, -1]$ satisfying
$f(-1)=0$ and for each $\eps>0$, there exists a $p\in \Pi$ of the form
$$p(x) =\sum_{k=0}^n a_k x^k$$ simultaneously satisfying
$$\|f-p\|_{[a, -1]} <\eps$$ and $$\sum_{k=0}^n |a_k| < \eps.$$

\sect{Some Multivariate Density Results}
In the section we consider applications of the results of the
previous sections to multivariate functions.

\medskip\noindent
{\bf Example 6.1.} We start with an application of the
Stone-Weierstrass theorem. Let $h_1\nek h_m$ be any $m$ fixed
real-valued continuous functions defined on $X$, a compact set.
Let $$\calM =\span\{g\big(\sum_{i=1}^m a_ih_i\big):
\bfa=(a_1\nek a_m)\in \RR^m,\, g\in C(\RR)\}.$$ When is $\calM$
dense in $C(X)$?

\proclaim Proposition \label{6.1}. $\overline{\calM}=C(X)$ if and
only if for each $x, y\in X$, $ x\ne  y$, there exists $i\in
\{1\nek m\}$ such that $h_i(x)\ne h_i(y)$.

\pf If there exists an $x\ne y$ for which $h_i(x) =
h_i(y)$, $i=1,...,m$, then for every $f\in \calM$ we have
$f(x)=f(y)$ and obviously $\overline{\calM} \ne C(X)$. On the
other hand, assume that for each $x, y\in X$, $ x\ne  y$, there
exists $i\in \{1\nek m\}$ such that $h_i(x)\ne h_i(y)$. Consider
the linear span of the set $$\big(\sum_{i=1}^m a_ih_i\big)^k$$ as
we vary over all $\bfa\in \RR^m$ and $k=0,1,2,\ldots$. This is an
algebra generated by $$h_1^{\ell_1}\cdots h_m^{\ell_m}$$ where the
$\ell_i$ are non-negative integers. Furthermore this algebra
contains the constant function and separates points. Thus the
density follows from the Stone--Weierstrass theorem.\eop

\noindent {\bf Example 6.2.} Is it true that for arbitrary compact
sets $X$ and $Y$ we always have that $\span \{C(X) {\times}
C(Y)\}$ is dense in $C(X {\times} Y)$? If $X$ and $Y$ are compact
subsets of $\RR$ this follows from the fact that algebraic
polynomials are dense in $C(X {\times} Y)$, and each algebraic
polynomial is a linear combination of products of monomials in $X$
and monomials in $Y$. Similarly to Example 6.1 we have:

\proclaim Dieudonn\'e Theorem \label{6.2}. If $X$ and $Y$ are
compact, then the linear space $\span\{ C(X){\times} C(Y)\}$ is dense in
$C(X{\times} Y)$.

\pf  For $f\in C(X)$ and $g\in C(Y)$ the function
$(f{\times} g)(x,y) = f(x)g(y)$ is in $C(X{\times} Y)$.
Furthermore all finite sums of the form $f_1{\times} g_1 +\cdots +
f_m{\times} g_m$ clearly form a subalgebra of $C(X{\times} Y)$
that contains the constant function and separates points. Thus by
the Stone-Weierstrass theorem $\span\{ C(X) {\times} C(Y)\}$  is
dense in $C(X{\times} Y)$. \eop

This theorem, originally proven in Dieudonn\'e [1937] by other
methods, easily extends to a product of any finite number of
compact spaces. It may also be found in Nachbin [1967] and Prolla
[1977].

\bigskip\noindent
{\bf Example 6.3.} Here is a simple application of the functional
analytic approach to density. We consider $C(\RR^n)$ with the
topology of uniform convergence on compacta. As we recall, the set
of functions for which the span of all their translates are not
dense in $C(\RR^n)$ are called {\it mean-periodic} functions.
There is no known characterization of mean-periodic functions in
$C(\RR^n)$ for $n\ge 2$. However not many functions can be
mean-periodic. For example

\proclaim Proposition \label{6.3}. If $g\in C(\RR^n)\cap
L^1(\RR^n)$, $g\ne 0$, then
$$C(\RR^n)=\overline{\span}\{g(\cdot-\bfa):\bfa\in \RR^n\}.$$

\pf The continuous linear functionals on $C(\RR^n)$ are
represented by Borel measures of bounded total variation and
compact support. If the above space is not dense in $C(\RR^n)$,
then there exists such a nontrivial measure $\mu$ satisfying
$$\int_{\RR^n} g(\bfx-\bfa)\dd \mu(\bfx)=0,$$ for all
$\bfa\in\RR^n$. Both $g$ and $\mu$ have ``nice'' Fourier
transforms. Since the above is a convolution we must have
$$\hatg(\bfw)\hatmu(\bfw)=0.$$ Now $\hatmu$ is an entire
function, while $\hatg$ is continuous. Since $\hatg$ must vanish
where $\hatmu\ne 0$, it follows that $\hatg=0$ and thus $g=0$, a
contradiction. \eop

\noindent {\bf Example 6.4.} Let $\langle\cdot\,,\cdot\rangle$
denote the usual inner (scalar) product on $\RR^n$. Applying
Propositions \recall{6.1} and \recall{5.8} we prove the following
result.

\proclaim Proposition \label{6.4}. For each $\sig\in C(\RR)$ $$
\span\{\sigma(\langle\bfa\,,\cdot\rangle +b)\,:\, \bfa\in \RR^n,
b\in \RR\}$$ is dense in $C(\RR^n)$ (uniform convergence on
compacta) if and only if $\sig$ is not a polynomial.

\pf If $\sig$ is a polynomial of degree $m$, then each
$\sigma(\langle\bfa\,, \cdot\rangle +b)$ is contained in the space
of polynomials of total degree at most $m$ on $\RR^n$, and thus the
above span is certainly not dense in $C(\RR^n)$.

Assume $\sig$ is not a polynomial. Choose an $f$ in $C(\RR^n)$,
$X$ any compact subset of $\RR^n$, and $\eps>0$. From an
application of Proposition \recall{6.1} we have the existence of
$g_k\in C(\RR)$ and $\bfa^k\in \RR^n$, $k=1\nek m$, such that
$$\left|f(\bfx) -\sum_{k=1}^m
g_k(\langle\bfa^k\,,\bfx\rangle)\right|<\eps$$ for all $\bfx\in
X$. Let $[c,d]$ be a finite interval of $\RR$ containing all
values $\langle\bfa^k\,,\bfx\rangle$ for $\bfx\in X$ and $k=1\nek
m$, i.e., $$\union_{k=1}^m \{ \langle\bfa^k\,,
\bfx\rangle\,:\,\bfx\in X\}\subseteq [c,d].$$ From Proposition
\recall{5.8} we have the existence of $c_{ik}, \alp_{ik},
\beta_{ik}\in \RR$, $i=1\nek n_k$, $k=1\nek m$, for which $$\left|
g_k(t) -\sum_{i=1}^{n_k} c_{ik} \sig(\alp_{ik}t +
\beta_{ik})\right| < {\eps\over m}$$ for all $t\in [c,d]$ and
$k=1\nek m$. Thus for all $\bfx\in X$ $$\left|f(\bfx)
-\sum_{k=1}^m \sum_{i=1}^{n_k}
c_{ik}\sig(\alp_{ik}\langle\bfa^k\,, \bfx\rangle
+\beta_{ik})\right|<2\eps,$$ which proves the density. \eop

Proposition \recall{6.4} is a basic result in one of the models of
neural network theory, see Leshno, Lin, Pinkus, Schocken [1993]
and Pinkus [1999].

\bigskip\noindent
{\bf Example 6.5.} {\it Ridge Functions.} Ridge functions were
considered in the previous example. They are functions of the form
$g(\langle\bfa\,, \cdot\rangle)$ for some fixed `direction' $\bfa$
and some function $g\in C(\RR)$. They are functions constant on
the hyperplanes $\{\langle\bfa\,, \bfx\rangle =t\}$ for every
$t\in \RR$.

Let $\Ome$ be a subset of $\RR^n$. In what follows we assume that
$\Ome$ is a subset of $S^{n-1}$, i.e., all elements of $\Ome$ are
of norm 1. (This is simply a convenient normalization.) The
question we ask is: What are necessary and sufficient conditions on
$\Ome$ such that the set of all ridge functions with directions
from $\Ome$ are dense in $C(\RR^n)$. The result we prove is due to
Vostrecov, Kreines [1961], see also Lin, Pinkus [1993]. We will
apply both the Weierstrass theorem and the Riesz representation
theorem in obtaining these conditions.

Let $$\calM(\Ome)= \span\{g(\langle\bfa\,, \cdot\rangle):\,\bfa\in
\Ome,\, g\in C(\RR)\}.$$ Note that we vary over all $\bfa\in
\Ome$ and {\it all} $g\in C(\RR)$.

\proclaim Theorem \label{6.5}. The linear space $\calM(\Ome)$ is
dense in $C(\RR^n)$ in the topology of uniform convergence on
compacta if and only if the only homogeneous polynomial (of $n$
variables) that vanishes identically on $\Ome$ is the zero
polynomial.

\pf ($\Rightarrow$). Assume there exists a nontrivial
homogeneous polynomial $p$ of degree $k$ that vanishes on $\Ome$.
Let $$p(\bfy) =\sum_{|\bfm|=k} b_{\bfm} \bfy^\bfm,$$ where
$\bfm=(m_1\nek m_n)\in \ZZ^n_+$, $|\bfm| = m_1+\cdots + m_n$ and
$\bfy^\bfm = y_1^{m_1}\cdots y_n^{m_n}$.

Choose any $\phi\in C^\infty_0(\RR^n)$, $\phi\ne 0$. For each
$\bfm\in\ZZ^n_+$, $|\bfm|=k$, set $$D^\bfm={\partial^k\over
\partial x_1^{m_1}\cdots \partial x_n^{m_n}},$$ and define
$$\psi(\bfx)=\sum_{|\bfm|=k}b_{\bfm}D^{\bfm}\phi(\bfx).$$ Note
that $\psi\in C^\infty_0(\RR^n)$, $\psi\ne 0$, ($\supp \psi\incl
\supp\phi$), and $$\hatpsi=i^k\hatphi p$$ where
$\,\widehat{\cdot}\,$ denotes the Fourier transform. As $p$ is
homogeneous, $p(\lambda \bfa)= \hatpsi(\lambda \bfa)=0$ for all
$\bfa\in \Ome$ and $\lambda\in \RR$.

We claim that $$\int_{\RR^n}g(\langle\bfa\,,
\bfx\rangle)\psi(\bfx)\dd \bfx=0$$ for all $\bfa\in \Ome$ and $g\in
C(\RR)$, i.e., the nontrivial linear functional defined by
integrating against $\psi$ annihilates $\calM(\Ome)$. From the
Riesz representation theorem this implies that $\calM(\Ome)$ is
not dense in $C(\RR^n)$.

We prove this as follows. For  $\bfa\in \Ome$ we write $$0=\hatpsi
(\lambda \bfa) = {1\over (2\pi)^{n/2}}\int_{\RR^n} \psi(\bfx)
e^{-i\lambda \langle\bfa\,,\bfx\rangle}\dd \bfx$$ $$={1\over
(2\pi)^{n/2}}\int_{-\infty}^\infty
\Big[\int_{\langle\bfa\,,\bfx\rangle=t}\psi(\bfx)d\bfx \Big]
e^{-i\lambda t}\dd t. $$ Since this holds for all $\lambda\in
\RR$, we have that $$\int_{\langle\bfa\,,
\bfx\rangle=t}\psi(\bfx)\dd \bfx=0$$ for all $t$. Thus for any $g\in
C(\RR)$,
$$\int_{\RR^n}g(\langle\bfa\,,\bfx\rangle)\psi(\bfx)\dd \bfx=
\int_{-\infty}^\infty \Big[\int_{\langle\bfa\,, \bfx\rangle=t}
\psi(\bfx)d\bfx\Big] g(t)\dd t=0.$$

\medskip \noindent ($\Leftarrow$). Assume that for a given $k\in \NN$
no nontrivial homogeneous polynomial $p$ of degree $k$ vanishes
identically on $\Ome$. We will prove that $\calM(\Ome)$ includes
all homogeneous polynomials of degree $k$ (and thus all
polynomials of degree at most $k$). If the above holds for all
$k\in \ZZ_+$ it then follows that $\calM(\Ome)$ contains all
polynomials and therefore $\overline{\calM(\Ome)}=C(\RR^n)$.

Let $\bfa\in\Ome$ and set $g(\langle\bfa\,, \bfx\rangle)=
(\langle\bfa\,,\bfx\rangle)^k$, whence
$(\langle\bfa\,,\bfx\rangle)^k\in \calM(\Ome)$. Since
$D^{{\bfm}_1}\bfx^{{\bfm}_2}=\del_{{\bfm}_1,{\bfm}_2} k!$, for
${\bfm}_1,{\bfm}_2\in \ZZ^n_+$, $|{\bfm}_1|=|{\bfm}_2|=k$, it
easily follows that every linear functional $\ell$ on the finite
dimensional linear space of homogeneous polynomials $H^n_k$ of
degree $k$ may be represented by some $q\in H^n_k$ via
$$\ell(p)=q(D)p$$ for each $p\in H^n_k$.

For any given $q\in H^n_k$,
$$q(D)\,(\langle\bfa\,,\bfx\rangle)^k=k!\,q(\bfa).$$ If the
linear functional $\ell$ annihilates $(\langle\bfa\,,
\bfx\rangle)^k$ for all $\bfa\in \Ome$, then its representor $q\in
H^n_k$ vanishes on $\Ome$. By assumption this implies that $q=0$.
The fact that no nontrivial linear functional on $H^n_k$
annihilates $(\langle\bfa\,,\bfx\rangle)^k$ for all $\bfa\in \Ome$
implies $$H^n_k=\span \{(\langle\bfa\,,\bfx\rangle)^k:\, \bfa\in
\Ome\}.$$ Thus $H^n_k\incl \calM(\Ome)$. \eop

\noindent {\bf Example 6.6.} There are other results of the same
general flavor as that found in Example 6.4. For example, assume
$\|\cdot\|$ is the usual Euclidean norm on $\RR^n$. Then we have
from Pinkus [1996] the following two results.

\proclaim Proposition \label{6.6}. For each $\sig\in C(\RR_+)$ $$
\span\{\sigma(\rho\|\cdot -\bfa\|)\,:\, \rho > 0,\, \bfa\in
\RR^n\}$$ is dense in $C(\RR^n)$ (uniform convergence on compacta)
if and only if $\sig$ is not an even polynomial.\nopf

\proclaim Proposition \label{6.7}. For each $\sig\in C(\RR)$
$$\span\{\sig\left(a \prod_{i=1}^n (\cdot-b_i)\right): \,a,b_1\nek
b_n\in \RR\}$$ is dense in $C(\RR^n)$ (uniform convergence on
compacta) if and only if $\sig$ is not of the form
$$\sig(t)=\sum_{j=0}^r c_{0j}t^j +\sum_{j=1}^r \sum_{i=1}^{n-1}
c_{ij} t^j (\ln |t|)^i$$ for some finite $r$ and coefficients
$(c_{ij})$.\nopf

\noindent {\bf Example 6.7.} A very interesting result, with
applications in Radon transform theory, is that conjectured by Lin
and Pinkus and proved by Agranovsky, Quinto [1996]. It
characterizes the set of centers of radial functions needed for
density. The complete answer is only known in $\RR^2$.

\proclaim Theorem \label{6.8}. Let $\calA\subseteq \RR^2$. Then
$$\span\{ g(\|\cdot-\bfa\|): \,g\in C(\RR),\, \bfa\in \calA\}$$
is not dense in $C(\RR^2)$ (uniform convergence on compacta) if
and only if $\calA$ is composed of a finite number of points together with a
subset of a set of straight lines having a common intersection
point and where the angles between each of the lines is a rational
multiple of $\pi$ (a Coxeter system of lines).\nopf

\noindent {\bf Example 6.8.} The functions $\|\cdot -\bfa\|^2$ are
shifts of the polynomial $q(\bfx) = \sum_{i=1}^n x_i^2$. Assume we
are given an arbitrary polynomial $p$. Under what exact conditions
do we have that
$$\span\{ g(p(\,\cdot - \bfa)):\, g\in C(\RR),\, \bfa\in \RR^n\}$$
is dense in $C(\RR^n)$? This next result, as well as variations
thereof, can be found in Pinkus, Wajnryb [1995].

\proclaim Theorem \label{6.9}. Let $p$ be an arbitrary polynomial
in $\RR^n$. Then for $n=1, 2, 3$,
$$\span\{ g(p(\,\cdot - \bfa)):\, g\in C(\RR),\, \bfa\in \RR^n\}$$
is dense in $C(\RR^n)$ (uniform convergence on compacta) if and
only if
$$\span\{ p(\,\cdot - \bfa):\, \bfa\in \RR^n\}$$
separates points.\nopf

By ``separates points'' we mean that for any given $\bfx, \bfy\in
\RR^n$, $\bfx \ne \bfy$, there exists a $\bfa\in \RR^n$ for which
$$p(\bfx-\bfa) \ne p(\bfy-\bfa).$$
This condition is obviously necessary. The sufficiency is far from
trivial. For $n=4$ it is also sufficient if $p$ is a homogeneous
polynomial. However for $n\ge 4$ this condition is not always
sufficient.

\bigskip \noindent
{\bf Example 6.9.} {\it M\"untz's Theorem.} The M\"untz problem in
the multivariate setting is significantly more difficult than in the
univariate setting. Some sufficient conditions have been given, but the problem
still remains very much open. The interested reader is urged to look at
Bloom [1992] and Kro\'o [1994] and references therein.

\medskip

\References

\refB Achieser, N.~I.; Theory of Approximation; Frederick Ungar
(New York); 1956; Originally published in Russian in 1947.

\refJ Agranovsky, M.~L., Quinto, E.~T.; Injectivity sets for the
Radon transform over circles and complete systems of radial
functions; J.\ Funct.\ Anal.; 139; 1996; 383--414;

\refJ Ahiezer, N.~I.; On the weighted approximation of continuous
functions by polynomials on the entire number axis; Uspehi Mat.\
Nauk; 11; 1956; 3--43; Also appears in translation in English in
{\sl American Math.\ Soc.\ Transl., Series 2} {\bf 22}, 95--137,
1962.

\refB Altomare, F., Campiti, M.; Korovkin-type Approximation
Theory and its Applications; Walter de Gruyter (Berlin - New
York); 1994;

\refJ Atzmon, A., Olevski\v{\i}, A.; Completeness of integer
translates in function spaces on $\RR$; J.\ Approx.\ Theory; 87;
1996; 291--327;

\refB Baillaud, B., Bourget, H.; Correspondance d'Hermite et de
Stieltjes, Tome II; Gauthier-Villars (Paris); 1905;

\refJ Banach, S.; Sur les fonctionelles lin\'eaires; Studia Math.;
1; 1929; 211--216 and 223--239;

\refB Banach, S.; Th\'eorie des Op\'erations Lin\'eaires; Hafner
(New York); 1932;

\refJ Bernstein, S.~N.; Sur l'ordre de la meilleure approximation
des functions continues par les polyn\^ome de degr\'e donn\'e;
Mem.\ Cl.\ Sci.\ Acad.\ Roy.\ Belg.; 4; 1912; 1--103;

\refJ Bernstein, S.~N.; D\'emonstration du th\'eor\`eme de
Weierstrass fond\'ee sur le calcul des
   probabilit\'es;
Comm.\ Soc.\ Math.\ Kharkow; 13; 1912/13; 1--2; Also appears in
Russian translation in Bernstein's Collected Works.

\refQ Bernstein, S.; Sur les recherches r\'ecentes relatives \`a
la meilleure approximation des
   fonctions continues par des polyn\^omes;
(Proceedings of the Fifth International Congress of
Mathematicians, Vol.~I), (Cambridge, 22-28 August 1912),
E.~W.~Hobson, A.~E.~H.~Love (eds.), Cambridge (England); 1913;
256--266; Also appears in Russian translation in Bernstein's
Collected Works.

\refJ Bernstein, S.; Le probl\`eme de l'approximation des
fonctions continues sur tout l'axe r\'eel
   et l'une de ses applications;
Bull.\ Soc.\ Math.\ France; 52; 1924; 399--410;

\refJ Bloom, T.; A multivariable version of the M\"untz-Sz\'asz theorem.
The Madison Symposium on Complex Analysis, (Madison, WI, 1991) ;
Contemp. Math.; 137; 1992; 85--92;

\refJ Bohman, H.; On approximation of continuous and of analytic
functions; Ark.\ Mat.; 2; 1952; 43--56;

\refJ Bonsall, F.~F.; Linear operators in complete positive cones;
Proc.\ London Math.\ Soc.; 8; 1958; 53--75;

\refJ Boor, C. de; On uniform approximation by splines; J.\
Approx.\ Theory; 1; 1968; 219--235;

\refJ Borwein, P.; Zeros of Chebyshev polynomials in Markov
systems; J.\ Approx.\ Theory; 63; 1990; 56--64;

\refB Borwein, P., Erd\'elyi, T.; Polynomials and Polynomial
Inequalities; Springer-Verlag (New York); 1995;

\refJ Borwein, P., Erd\'elyi, T.; Dense Markov spaces and
unbounded Bernstein inequalities; J.\ Approx.\ Theory; 81; 1995;
66--77;

\refB Buck, R.~C.; Studies in Modern Analysis, Vol.~1;
Mathematical Association of America (Washington D.~C); 1962;
Contains a reprint of Stone [1948].

\refJ Burckel, R.~B., Saeki, S.; An elementary proof of the
M\"untz-Sz\'asz theorem; Expo.\ Math.; 4; 1983; 335--341;

\refB Cheney, E.~W.; Introduction to Approximation Theory;
McGraw-Hill (New York); 1966;

\refJ Conrey, J.~B.; The Riemann Hypothesis; Notices of the
A.~M.~S.; 50; 2003; 341--353;

\refJ Dieudonn\'e, J.; Sur les fonctions continues num\'erique
d\'efinies dans une produit de deux
   espaces compacts;
Comptes Rendus Acad.\ Sci.\ Paris; 205; 1937; 593--595;

\refB Feinerman, R.~P., Newman, D.~J.; Polynomial Approximation;
Wil\-liams and Wilkins (Baltimore); 1974;

\refJ Fej\'er, L.; Sur les fonctions born\'ees et int\`egrables;
Comptes Rendus Acad.\ Sci.\ Paris; 131; 1900; 984--987;

\refJ Golitschek, M. v.; A short proof of M\"untz's theorem; J.\
Approx.\ Theory; 39; 1983; 394--395;

\refJ Hahn, H.; \"Uber lineare Gleichungssysteme in linearen
R\"aumen; J.\ Reine Angew.\ Math.; 157; 1927; 214--229;

\refJ Helly, E.; \"Uber lineare Functionaloperationen; Sitzung.
der \"Oster.\ Akad.\ Wissen.\ (Wien); 121; 1912; 265--297;

\refJ Korovkin, P.~P.; On convergence of linear positive operators
in the space of continuous functions; Dokl.\ Akad.\ Nauk SSSR; 90;
1953; 961--964;

\refB Korovkin, P.~P.; Linear Operators and Approximation Theory;
Hindustan Publ. Corp. (Delhi); 1960; (The Russian original
appeared in 1959.)

\refJ Kro\'o, A.; A geometric approach to the multivariate M\"untz
problem; Proc. Amer. Math.Soc.; 121; 1994; 199--208;

\refJ Lebesgue, H.; Sur l'approximation des fonctions; Bull.\
Sciences Math.; 22; 1898; 278--287;

\refJ Lerch, M.; O hlavni vete theorie funkci vytvorujicich
   (On the main theorem on generating functions);
Rozpravy Ceske Akademie v.~Praze; 1; 1892; 681--685;

\refJ Lerch, M.; Sur un point de la th\`eorie des fonctions
g\'en\'eratices d'Abel; Acta Math.; 27; 1903; 339--351;

\refJ Leshno, M., Lin, V.~Ya., Pinkus, A., Schocken, S.;
Multilayer feedforward networks with a non-polynomial activation
function can
   approximate any function;
Neural Networks; 6; 1993; 861--867;

\refB Levin, B.~Ja.; Distribution of Zeros of Entire Functions;
Transl.\ Math.\ Monographs, {\bf 5}, Amer.\ Math.\ Soc.
(Providence); 1964;

\refB Levin, B.~Ja.; Lectures on Entire Functions; Transl.\ Math.\
Monographs, {\bf 150}, Amer.\ Math.\ Soc. (Providence); 1996;

\refB Levinson, N.; Gap and Density Theorems; Amer.\ Math.\ Soc.,
Colloquium Publ., {\bf 26} (New York); 1940;

\refJ Lin, V.~Ya., Pinkus, A.; Fundamentality of ridge functions;
J.\ Approx.\ Theory; 75; 1993; 295--311;

\refB Lorentz, G.~G., v.~Golitschek, M., Makovoz, Y.; Constructive
Approximation. Advanced problems; Springer Verlag (Berlin); 1996;

\refQ Luxemburg, W.~A.~J.; M\"untz-Sz\'asz type approximation
results and the Paley-Weiner theorem; (Approximation Theory II),
G.~G.~Lorentz, C.~K.~Chui, L.~L.~Schumaker (eds.), Academic Press
(New York); 1976; 437--448;

\refJ Luxemburg, W.~A.~J., Korevaar, J.; Entire functions and
M\"untz-Sz\'asz type approximation; Trans.\ Amer.\ Math.\ Soc.;
157; 1971; 23--37;

\refJ Mergelyan, S.~N.; Weighted approximation by polynomials;
Uspehi Mat.\ Nauk; 11; 1956; 107--152; Also appears in translation
in English in {\sl American Math.\ Soc.\ Transl., Series 2} {\bf
10}, 59--106, 1958.

\refB Mhaskar, H.~N.; Introduction to the Theory of Weighted
Polynomial Approximation; World Scientific (Singapore); 1996;

\refQ M\"untz, C.~H.; \"Uber den Approximationssatz von
Weierstrass; (Mathematische Abhandlungen Hermann Amandus Schwarz
zu seinem f\"unfzigj\"ahrigen Doktojubil\"aum am 6. August 1914
gewidmet von Freunden und Sch\"ulern), C. Carath\'eodory, G.
Hessenberg, et al. (eds.), Springer (Berlin); 1914; 303--312;

\refB Nachbin, L.; Elements of Approximation Theory; Van Nostrand
(Princeton); 1967;

\refJ Nikolski, N.; Remarks concerning completeness of translates
in function spaces; J.\ Approx.\ Theory; 98; 1999; 303--315;

\refJ Nussbaum, R.~D., Walsh, B.; Approximation by polynomials
with nonnegative coefficients and the spectral theory of positive
operators; Trans.\ Amer.\ Math.\ Soc.; 350; 1998; 2367--2391;

\refB Orlicz, W.; Linear Functional Analysis; World Scientific
(Singapore); 1992;

\refB Paley, R., Wiener, N.; Fourier Transforms in the Complex
Domain; Amer.\ Math.\ Soc., Colloquium Publ., {\bf 19} (New York);
1934;

\refJ Picard, E.; Sur la repr\'esentation approch\'ee des
fonctions; Comptes Rendus Acad.\ Sci.\ Paris; 112; 1891; 183--186;

\refJ Pinkus, A.; TDI-subspaces of $C(\RR^d)$ and some density
problems from neural networks; J.\ Approx.\ Theory; 85; 1996;
269--287;

\refJ Pinkus, A.; Approximation theory of the MLP model in neural
networks; Acta Numerica; 8; 1999; 143--195;

\refJ Pinkus, A.; Weierstrass and approximation theory; J.\
Approx.\ Theory; 107; 2000; 1--66;

\refJ Pinkus, A.; The Weierstrass approximation theorems; 
Surveys in Approx. Theory;  ; 2005; (to appear);

\refJ Pinkus, A., Wajnryb, B.; A Problem of Approximation Using
Multivariate Polynomials; Uspekhi Mat.\ Nauk; 50; 1995; 89--110; =
{\sl Russian Math.~Surveys} {\bf 50} (1995), 319--340.

\refB Prolla, J.~B.; Approximation of Vector Valued Functions;
North-Holland (Amsterdam); 1977;

\refB Prolla, J.~B.; Weierstrass--Stone, the Theorem; Peter Lang
(Frankfurt); 1993;

\refJ Riesz, F.; Sur les op\'erations fonctionnelles lin\'eaires;
Comptes Rendus Acad.\ Sci.\ Paris; 149; 1909; 974--977; Also
appears in {\sl \OE uvres}, Volume I, 400--402.

\refJ Riesz, F.; Sur certains syst\`emes d'equations fonctionelles
et l'approxim\-ation des
   fonctions continues;
Comptes Rendus Acad.\ Sci.\ Paris; 150; 1910; 674--677; Also
appears in {\sl \OE uvres} of F.~Riesz. This paper is on pages
403--406 of Volume I. Unfortunately, in some copies of the \OE
uvres, there was a mix-up with a previous paper, and the pages are
403 and 404 followed by 398 and 399.

\refJ Riesz, F.; Sur certains syst\`emes singuliers d'\'equations
int\'egrales; Ann.\ Sci.\ de l'Ecole Norm.\ Sup.; 28; 1911;
33--62; Also appears in {\sl \OE uvres}, Volume II, 798--827.

\refJ Rogers, L.~C.~G.; A simple proof of M\"untz's theorem;
Math.\ Proc.\ Camb.\ Phil.\ Soc.; 90; 1981; 1--3;

\refB Rudin, W.; Real and Complex Analysis; McGraw-Hill (New
York); 1966;

\refJ Runge, C.; \"Uber die Darstellung willk\"urlicher
Functionen; Acta Math.; 7; 1885/86; 387--392;

\refD Schmidt, E.; Entwicklung willk\"urlicher Funktionen nach
Systemen vor\-ge\-schrie\-be\-ner; G\"ottingen (Germany); 1905;

\refB Schwartz, L.; \'Etude des sommes d'exponentielles r\'eelles;
Hermann et Cie. (Paris); 1943;

\refJ Schwartz, L.; Sur certaines familles non fondamentales de
fonctions continues; Bull.\ Soc.\ Math.\ France; 72; 1944;
141--145;

\refJ Schwartz, L.; Th\'eorie g\'en\'erale des fonctions
moyenne-p\'eriodiques; Ann.\ Math.; 48; 1947; 857--928;

\refJ Stone, M.~H.; Applications of the theory of Boolean rings to
general topology; Trans.\ Amer.\ Math.\ Soc.; 41; 1937; 375--481;

\refJ Stone, M.~H.; A generalized Weierstrass approximation
theorem; Math.\ Magazine; 21; 1948; 167--184, 237--254;

\refJ Sz\'asz, O.; \"Uber die Approximation stetiger Funktionen
durch lineare Aggregate von Potenzen ; Math.\ Ann.; 77; 1916;
482--496;

\refJ Toland, J.~F.; Self-adjoint operators and cones; J.\ London
Math.\ Soc.; 53; 1996; 167--183;

\refJ Vostrecov, B.~A., Kreines, M.~A.; Approximation of
continuous functions by superpositions of plane waves;  Dokl.\
Akad.\ Nauk SSSR; 140; 1961; 1237--1240; = {\sl Soviet Math.
Dokl.} {\bf 2} (1961), 1326--1329.

\refJ Weierstrass, K.; \"Uber die analytische Darstellbarkeit
sogenannter will\-k\"ur\-li\-cher Functionen einer reellen
Ver\"anderlichen; Sitzungsberichte der Akademie zu Berlin; ; 1885;
633--639 and 789--805; (This appeared in two parts. An expanded
version of this paper with ten additional pages also appeared in
{\sl Mathematische Werke}, {\bf Vol.~3}, 1--37, Mayer \& M\"uller,
Berlin, 1903.)

\refJ Weierstrass, K.; Sur la possibilit\'e d'une repr\'esentation
analytique des fonctions dites arbitraires d'une variable
r\'eelle; J.\ Math.\ Pure et Appl.; 2; 1886; 105--113 and 115-138;
(This is a translation of Weierstrass [1885] and, as the original,
it appeared in two parts and in subsequent issues, but under the
same title. This journal was, at the time, called {\sl Journal de
Liouville})


{

\bigskip\obeylines
Allan Pinkus
Department of Mathematics
Technion
Haifa, 32000
Israel
{\tt pinkus@tx.technion.ac.il}
{\tt http://www.technion.ac.il/\~{}pinkus/}

}

\bye